\documentclass{amsart}
\usepackage{amsmath,amsfonts,amssymb,amsthm,epsfig,amscd,psfrag,latexsym,comment}
\usepackage[all]{xy}

\input xy
\xyoption{all}

\newcommand{\puteps}[2][0.5]
{\includegraphics[scale=#1]{#2.eps}}

\newtheorem{thm}{Theorem}[section]
\newtheorem{lem}[thm]{Lemma}

\newtheorem{cor}[thm]{Corollary}
\newtheorem{prop}[thm]{Proposition}

\newtheorem{prob}[thm]{Problem}

\newtheorem{Specialthm}{Theorem}

\theoremstyle{definition}
\newtheorem{Remark}[thm]{Remark}

\def\Pic{{\mbox{Pic}}}
\def\Gr{{\sf Gr}}
\def\sC{{\mathcal{C}}}
\def\sD{{\mathcal{D}}}
\def\C{{\mathbb C}}

\def\tS{\widetilde{S}}
\def\tT{\widetilde{T}}
\def\tp{\widetilde{p}}
\def\Z{{\mathbb Z}}
\def\Q{{\mathbb Q}}

\def\P{{\mathbb P}}
\def\H{{\mathcal H}}
\def\O{{\mathcal O}}
\def\Cu{{\mathcal C}}

\def\dim{{\mbox{dim}}}
\def\Aut{{\mbox{Aut}}}
\def\span{{\mbox{span}}}

\def\tX{{\widetilde{X}}}

\def\tU{{\widetilde{U}}}
\def\tV{{\widetilde{V}}}
\def\tY{{\widetilde{Y}}}

\newcounter{example}
\def\Example{
\addtocounter{example}{1}
{\bf Example \theexample. }}

\def\Spec{\mbox{Spec}}
\def\Spf{\mbox{Spf}}
\def\Sing{{\mbox{Sing}}}
\def\ip#1#2{\langle #1,#2 \rangle}
\def\sM{\mathcal{M}}
\def\oM{\overline{M}}
\def\oW{\overline{W}}
\def\soM{\overline{\mathcal{M}}}
\def\sA{{\mathcal A}}
\def\oX{{\overline{X}}}
\def\sX{{\mathcal{X}}}
\def\soX{{\overline{\mathcal{X}}}}
\def\oY{{\overline{Y}}}
\def\sY{{\mathcal{Y}}}

\def\ov{\overline{v}}

\begin{document}

\title[The Abelian Monodromy Extension Property for Families of Curves]{The Abelian Monodromy Extension Property for Families of Curves} 

\author{Sabin Cautis}
\email{scautis@math.harvard.edu}
\address{Department of Mathematics \\ Rice University \\ Houston, TX}

\begin{abstract}
Necessary and sufficient conditions are given (in terms of monodromy) for extending a family of smooth curves over an open subset $U \subset S$ to a family of stable curves over $S$. More precisely, we introduce the abelian monodromy extension (AME) property and show that the standard Deligne-Mumford compactification is the unique, maximal AME compactification of the moduli space of curves. We also show that the Baily-Borel compactification is the unique, maximal projective AME compactification of the moduli space of abelian varieties.
\end{abstract}

\date{\today}
\maketitle
\tableofcontents

\section{Introduction}\label{s0}

We work over $\C$. Fix a separated Deligne-Mumford (DM) stack $\sX$ and a compactification $\soX$. Given an open subset $U \subset S$ of a normal variety $S$ together with a regular map $U \rightarrow \sX$, when does this map extend to a regular map $S \rightarrow \soX$?

Let's illustrate why this can be interesting geometrically. We will denote by $\sM_{g,n}$ the moduli stack of $n$-pointed genus $g$ curves and by $M_{g,n}$ its coarse moduli space. Also, $\soM_{g,n}$ and $\oM_{g,n}$ will denote their respective Deligne-Mumford compactifications. 

Now take $\sX = \sM_g$ and $\soX = \soM_g$ (no marked points). The simplest case is when the dimension of $S$ is one. Since extending a map is a local problem we take $S=C$ to be a curve and $U = C \setminus \{p\}$ the complement of a point. Section 2 of \cite{dm} shows that a family of smooth curves over $C \setminus \{p\}$ extends to a family of stable curves over $C$ if and only if the associated Jacobian family extends to a family of semi-abelian varieties. On the other hand, sections 3.5--3.8 of \cite{gr} show that a family of abelian varieties over $C \setminus \{p\}$ extends to a family of semi-abelian varieties over $C$ if and only if the associated monodromy on the $H_1$ homology of a fibre in a small (analytic) neighbourhood of $p$ is unipotent. Combining these results one obtains:

\begin{thm}[Deligne-Mumford-Grothendieck]\label{thm:DMG} A family of smooth curves over $C \setminus \{p\}$ extends to a family of stable curves over $C$ if and only if the induced monodromy on $H_1$ of the fibres around an (analytic) neighbourhood of $p$ is unipotent. 
\end{thm}

\begin{figure}[ht]
\begin{center}
 \psfrag{D}{\footnotesize{$D$}}
 \psfrag{E1}{\footnotesize{$E_1$}}
 \psfrag{E2}{\footnotesize{$E_2$}}
 \psfrag{E3}{\footnotesize{$E_3$}}
 \psfrag{j=0}{\footnotesize{$j=0$}}
 \psfrag{j=1728}{\footnotesize{$j=1728$}}
 \psfrag{j=8}{\footnotesize{$j=\infty$}}
 \hspace{0 cm}\puteps[0.4]{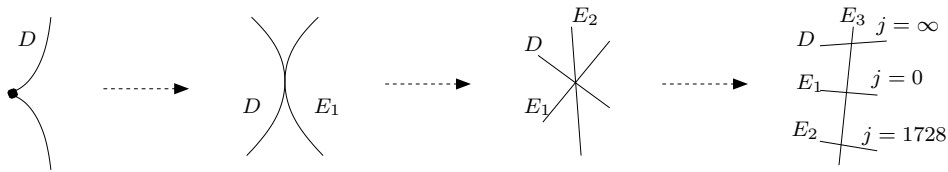}
\end{center}
\caption{$D = S \setminus U$ is the discrimant curve $4a^3+27b^2=0$.
\label{f1C}}
 \end{figure}

If $\dim(S) > 1$ you may need to blow up the base before you can extend $U \rightarrow \sM_g$. A neat example is the family $y^2=x^3+ax+b$ of elliptic curves over $a,b \in \C^2$. The fibres are smooth over $U = \C^2 \setminus \{(a,b): 4a^3+27b^2=0\}$ and stable over $\C^2 \setminus (0,0)$. The induced map $\C^2 \setminus (0,0) \rightarrow \oM_{1,1}$ does not extend over $(0,0)$. To extend it one has to blow up three times to get the surface $S'$ shown in Figure \ref{f1C}. The map $U \rightarrow \sM_{1,1}$ now extends to a morphism $S' \rightarrow \oM_{1,1}$ to the coarse moduli compactification which collapses the exceptional divisors $E_1$ and $E_2$ to the points with $j$ invariants $0$ and $1728$ and maps $E_3$ one-to-one onto $\oM_{1,1}$.

As this example illustrates, it was necessary to resolve $D = S \setminus U$ to a normal crossing divisor before the map extended. This is not quite a coincidence. In \cite{jo} it was shown that (see Remark \ref{rem:jo}):

\begin{thm}\label{thm:jo} Let $D = S \setminus U$ be a normal crossing divisor at $p$. Then a morphism $U \rightarrow \sM_{g,n}$ extends to a regular map $S \rightarrow \oM_{g,n}$ in a Zariski neighbourhood of $p$. 
\end{thm}

Theorem \ref{thm:jo} holds in arbitrary characteristic (which is significant) but if one works over $\C$ we have the following generalization (Corollary \ref{cor7}).

\begin{Specialthm}\label{specialthm1} 
Let $U \subset S$ be an open subvariety of an irreducible, normal variety $S$. A morphism $U \rightarrow \sM_{g,n}$ extends to a regular map $S \rightarrow \oM_{g,n}$ in a Zariski neighbourhood of $p \in S \setminus U$ if and only if the local monodromy around $p$ is virtually abelian. 
\end{Specialthm}

Notice that since the local fundamental group of the complement of a normal crossing divisor is abelian, Theorem \ref{specialthm1} implies \ref{thm:jo} over $\C$. 

\begin{Remark} The result in \cite{jo} actually gives an extension to $\soM_{g,n}$ (not just $\oM_{g,n}$) but only under the added condition that the morphism extends to a regular map $S \rightarrow \soM_{g,n}$ over the generic points of $D$. This condition is sufficient to lift the extension from $\oM_{g,n}$ to $\soM_{g,n}$ (see Corollary \ref{cor4}). 
\end{Remark}

One can ask if such extension results exist for other spaces. Another example involves the moduli stack of $g$-dimensional principally polarized abelian varieties $\sA_g$ and its Baily-Borel compactification $A_g^{BB}$. It was shown in \cite{b} that 

\begin{thm}\label{thm:b} 
Let $D = S \setminus U$ be a normal crossing divisor at $p$. Then a morphism $U \rightarrow \sA_g$ extends to a regular map $S \rightarrow A_g^{BB}$ in a Zariski neighbourhood of $p$. 
\end{thm}

As an immediate consequence of Theorem \ref{thm5} it follows that:

\begin{Specialthm}\label{specialthm2} 
Let $U \subset S$ be an open subvariety of an irreducible, normal variety $S$. A morphism $U \rightarrow \sA_g$ extends to a regular map $S \rightarrow A_g^{BB}$ in a Zariski neighbourhood of $p \in S \setminus U$ if the local monodromy around $p$ is virtually abelian.
\end{Specialthm}

Notice that in Thorem \ref{specialthm2} there is no ``only if'' part. This is because we don't know whether the local fundamental group of $\sA_g$ around a point in the boundary of $A_g^{BB}$ is virtually abelian. Understanding this local fundamental group is an interesting question in itself. 

Inspired by these results we introduce the abelian monodromy extension (AME) property for a pair $(\sX, \oX)$ consisting of a Deligne-Mumford stack $\sX$ and a compactification $\oX$ of its coarse scheme (section \ref{s1}). Roughly, we say $\oX$ is an AME compactification of $\sX$ if $U \rightarrow \sX$ extends to a regular map $S \rightarrow \oX$ whenever the image of the induced map $\pi_1(U) \rightarrow \pi_1(\sX)$ is virtually abelian. Here we think of $U$ and $S$ as being very small analytic neighbourhoods of a point (so this is a local condition on the domain but a global condition on the target $(\sX, \oX)$). Section \ref{s2} describes some basic properties of such AME compactifications and shows that if the AME compactification of $\sX$ exists then there is a unique maximal one which we denote $X_{ame}$ (Corollary \ref{cor6}).

Section \ref{s3} Theorem \ref{thm6} shows that $\oM_{g,n}$ is the maximal AME compactification of $\sM_{g,n}$ and derives Theorem \ref{specialthm1} as a corollary. Finally, section \ref{s4} Theorem \ref{thm5} proves that $A_g^{BB}$ is the maximal AME compactification of $\sA_g$ (among all projective ones).

\subsection{Future work}

In a future paper we will develop the theory of AME compactifications and give further geometric examples. We'll also explore a connection between AME compactifications and log canonical models which, among other things, gives a more direct proof of the maximality statement in Theorem \ref{thm5}.

\subsection{Acknowledgements}

I would like to thank Samuel Grushevsky, Joe Harris, Brendan Hassett, Klaus Hulek and Sean Keel for many insightful discussions. Johan de Jong and Jason Starr also made helpful suggestions early on while Richard Hain and J\'anos Koll\'ar were very kind in pointing out the utility of \cite{del} in proving Theorem \ref{thm5}. I appreciate the support and hospitality of the mathematics departments at Harvard and Rice University and of the Mittag-Leffler Institute where I was a visitor in the spring of 2007. I would also like to thank the referees for helpful, detailed suggestions in particular with respect to defining local monodromy and simplifying some proofs. 

\section{The Abelian Monodromy Extension (AME) Property}\label{s1}

In this section we define what it means for a variety (or stack) to have the abelian monodromy extension (AME) property. We also make a few remarks about the definition and give an example. 

All schemes are Noetherian. By a variety we mean an integral, separated (but possibly singular) scheme of finite type defined over $\C$. By morphism we will always mean a regular morphism. A complete scheme is a scheme proper over $\Spec(\C)$ and a compactification of a variety $X$ is a complete variety $\oX$ containing $X$ as an open dense subset. 

\subsection{Local and global monodromy}\label{s2.1}

Let $X$ be an open subvariety of a normal variety $\oX$. Next, consider an open subvariety $U \subset S$ of a normal variety $S$ together with a morphism $U \rightarrow X$. Fix a connected, reduced, complete subscheme $T \subset S$. We can always restrict to a smaller open subset $U$ if necessary so that $U$ is disjoint from $T$ as well as the singular locus of $S$. 

Let $D = S \setminus U$. Since we can find a blowup $\pi: S' \rightarrow S$ so that $\pi^{-1}(D \cup T)$, $\pi^{-1}(D)$ and $\pi^{-1}(T)$ are all simple, normal crossing divisors in $S'$ we can assume that $D \cup T$, $D$ and $T$ are all simple, normal crossing divisors. This means that around any point $p \in T$ there exist local coordinates such that $T$ is given by $x_1 \cdots x_i = 0$ and $D$ is given by $x_{i+1} \cdots x_{i+j} = 0$. Notice that we can choose $U \subset S$ so that this blowup leaves $U$ unchanged. 

Now choose a small analytic neighbourhood $V$ of $T$. One way to do this is to put a Riemannian metric $\rho$ on $S$ and take
$$V = V^{\rho}_\epsilon := \{p \in S: d_\rho(p, T) < \epsilon \}.$$
Then we can define the {\em local monodromy around $T$} as the image 
$$\mbox{Im} \left( \pi_1(V \cap U) \rightarrow \pi_1(X) \right).$$ 
Since $T$ is connected $V \cap U$ is connected and so the image of $\pi_1(V \cap U)$ is defined (up to conjugation) without having to choose a base point. 

There are three things to check to make sure local monodromy is well defined. The first is that local monodromy does not depend on the choice of $V$. To see this note we identify $D \cup T \subset S$ as sitting inside its normal bundle (which is a line bundle). This means that $V^\rho_\epsilon \cap (U \setminus T)$ deformation retracts onto $V^\rho_{\epsilon'} \cap (U \setminus T)$ for any $\epsilon' < \epsilon$ sufficiently small. Thus $\pi_1(V_\epsilon^\rho \cap (U \setminus T))$ and $\pi_1(V_{\epsilon'}^\rho \cap (U \setminus T))$ have isomorphic images in $\pi_1(X)$. Moreover, since for any other metric $\rho'$ we have $V^{\rho'}_{\epsilon'} \subset V^{\rho}_{\epsilon}$ for $\epsilon' \ll \epsilon$ we see that this image does not depend on the metric. 

Now the images of 
\begin{equation*}
\pi_1(V_\epsilon^\rho \cap (U \setminus T)) \rightarrow \pi_1(X) \text{ and }
\pi_1(V_\epsilon^\rho \cap U) \rightarrow \pi_1(X)
\end{equation*}
are the same since the map
$$\pi_1(V_\epsilon^\rho \cap (U \setminus T)) \rightarrow \pi_1(V_\epsilon^\rho \cap U)$$
induced by inclusion is surjective. Thus $\pi_1(V_{\epsilon}^{\rho} \cap U)$ and $\pi_1(V_{\epsilon'}^{\rho'} \cap U)$ have the same image for sufficiently small $\epsilon, \epsilon'$. 

The second thing is that local monodromy does not depend on the choice of $U$. This follows since if $U' \subset U$ is (in the analytic topology) a dense, open subset then the map $\pi_1(U') \rightarrow \pi_1(U)$ induced by the inclusion is surjective. Similarly, $\pi_1(V \cap U') \rightarrow \pi_1(V \cap U)$ is surjective and hence the images of 
\begin{equation*}
\pi_1(V \cap U) \rightarrow \pi_1(X) \text{ and } \pi_1(V \cap U') \rightarrow \pi_1(X)
\end{equation*}
are the same. 

The third thing to check is that local monodromy does not depend on the choice of the blowup $S'$. To see this consider a blowup $\pi: S' \rightarrow S$ so that $\pi^{-1}(D \cup T)$, $\pi^{-1}(D)$ and $\pi^{-1}(T)$ are also simple, normal crossing divisors. If you take a small neighbourhood $V$ of $T$ then $\pi^{-1}(V)$ is also a small neighbourhood of $\pi^{-1}(T)$ and then taking $U \subset S$ so that it avoids the blowup locus we get $\pi^{-1}(U) = U$ so that $V \cap U = \pi^{-1}(V) \cap \pi^{-1}(U)$ and hence the images of
\begin{equation*}
\pi_1(V \cap U) \rightarrow \pi_1(X) \text{ and } \pi_1(\pi^{-1}(V) \cap \pi^{-1}(U)) \rightarrow \pi_1(X)
\end{equation*}
are the same. 

Most commonly we will take $T$ to be a point $p \in S \setminus U$ to obtain the local monodromy around $p$. On the other hand, if $S$ is complete, we can take $T$ to be all of $S$ to obtain the {\em global monodromy of $S$} as the image of the induced map $\pi_1(U) \rightarrow \pi_1(X)$. 

The reason for this slightly strange terminology is that $X$ will often be a moduli space, for instance $\sM_g$ or $\sA_g$, so that the image of $\pi_1(U) \rightarrow \pi_1(X)$ can be identified with the usual notion of monodromy. For instance, $\pi_1(\sM_g)$ (the orbifold fundamental group) is the mapping class group $\Gamma_g$ so the image of $\pi_1(U) \rightarrow \pi_1(\sM_g)$ is the monodromy on $\pi_1$ of a fibre in the family of curves over $U$ corresponding to the morphism $U \rightarrow \sM_g$. Similarly, $\pi_1(\sA_g) = Sp_{2g}(\Z)$ so the image of $\pi_1(U) \rightarrow \pi_1(\sA_g)$ induced by a morphism $U \rightarrow \sA_g$ is the monodromy on $H_1$ of a fibre in the associated family of abelian varieties. 

\begin{Remark}
It does not actually matter what dense open $U$ we pick in the above definition of monodromy (as long as the map $U \rightarrow X$ is regular). This is because the map $\pi_1(V \cap U) \rightarrow \pi_1(V \cap U')$ is surjective if $U \subset U'$ is dense in the analytic topology so that the images of
\begin{equation*}
\pi_1(V \cap U) \rightarrow \pi_1(X) \text{ and } \pi_1(V \cap U') \rightarrow \pi_1(X)
\end{equation*}
are the same. Thus sometimes we will take $U$ to be the largest open subset where the map $S \dashrightarrow \oX$ is regular. 
\end{Remark}

\subsection{The AME property}

Given an open embedding of normal varieties $X \subset \oX$, the pair $(X,\oX)$ has {\em the abelian monodromy extension (AME) property} if given any $U \subset S$ as above the morphism $U \rightarrow X$ extends to a regular map $S \rightarrow \oX$ in a neighbourhood of $p$ whenever the local monodromy around $p$ is virtually abelian (a group is virtually abelian if it contains an abelian subgroup of finite index). In this case $\oX$ is complete (see Lemma \ref{lem1} below) and we say that $\oX$ is an {\em AME compactification} of $X$. We say $X$ has the AME property if it has an AME compactification. 

\begin{lem}\label{lem1} If the pair $(X,\oX)$ has the AME property then $\oX$ is complete.
\end{lem}
\begin{proof}
Let $\oY$ be a compactification of $\oX$ (any separated scheme of finite type over $\C$ has a compactification). Consider a blowup $\pi: \oY' \rightarrow \oY$ so that the boundary $B := \oY' \setminus \pi^{-1}(X)$ is a normal crossing divisor. Then the fundamental group of $V \cap \pi^{-1}(X)$ for a sufficiently small neighbourhood $V$ of any point $p \in B$ is abelian. 
Since $(X, \oX)$ has the AME property this means that the morphism $\pi^{-1}(X) \rightarrow X$ extends to a morphism $f: \oY' \rightarrow \oX$. Then $\oY'$ is complete and $f$ is dominant so $\oX$ must also be complete (since the image of a complete variety is complete). 
\end{proof}

\begin{Remark} Alternatively, one might want to define the AME property by requiring that the morphism extend whenever the local monodromy is abelian instead of virtually abelian. It turns out these two definitions are the same. To see this suppose the local monodromy is virtually abelian. The preimage under $\pi_1(U) \rightarrow \pi_1(X)$ of a finite index abelian subgroup inside $\mbox{Im}(\pi_1(U))$ is a finite index subgroup inside $\pi_1(U)$ and contains a normal finite index subgroup $H \subset \pi_1(U)$. Take the unramified cover $\pi: \widetilde{U} \rightarrow U$ corresponding to $H$ and let $\widetilde{S}$ be the normal closure of $S$ in the function field of $\widetilde{U}$ to get a finite morphism $\pi: \widetilde{S} \rightarrow S$. Then the local monodromy around any point of $\pi^{-1}(p)$ induced by $\widetilde{U} \rightarrow X$ is abelian and we get a regular extension $\widetilde{S} \rightarrow \oX$. By Lemma \ref{lem4} below this means we also get a regular extension $S \rightarrow \oX$. 

\begin{lem}\label{lem4} Consider the composition 
\[
\begin{xy}
\xymatrix{
\tS \ar[d]_{\pi} & \\
S \ar@{-->}[r]^h & X 
}
\end{xy}
\]
where $\tS$, $S$ are quasi-projective varieties, $S$ is normal and $X$ is complete. If $\pi$ is a finite, surjective morphism then the rational map $h$ is regular if and only if $h \circ \pi$ is regular. 
\begin{Remark}
Notice that this result does not hold if $X$ is a stack. 
\end{Remark}
\end{lem}
\begin{proof}
We need to show $h$ is regular assuming $h \circ \pi$ is regular. 

Let $S'$ be the closure of the image of $S$ in $S \times X$. Now if $U$ is the domain where $h$ is regular then $\pi^{-1}(U)$ is open and dense in $\tS$. Then the image of $(\pi, h \circ \pi): \tS \rightarrow S \times X$ is the closure of the image of $U$ which is $S'$. 

Thus $f: S' \subset S \times X \rightarrow S$ is finite (since $\pi: \tS \rightarrow S$ is finite) and birational (since $S'$ is the closure of the image of $S$). Since $S$ is normal this means $f$ is an isomorphism and thus $h$ is regular.  
\end{proof}
\end{Remark}

\Example Consider a smooth, Deligne-Mumford stable $n$-pointed curve $(C, \{p_1, \dots, p_n\})$. Then $C$ is an AME compactification of $C \setminus \cup_i p_i$. To see this the key observation is that $(C, \{p_1, \dots, p_n\})$ is stable precisely when $\pi_1(C \setminus \cup_i p_i)$ is not virtually abelian. 

More precisely, let $U \subset S$ be as above and suppose we have a morphism $U \rightarrow C \setminus \cup_i p_i$ which does not extend to a regular morphism around some point $p \in S \setminus U$. We want to show that this implies the local monodromy around $p$ is not virtually abelian. 

First resolve the map $S \dashrightarrow C$ by some blowup $\pi: S' \rightarrow S$ to get a morphism $f: S' \rightarrow C$. Since the morphism $U \rightarrow C \setminus \cup_i p_i$ did not extend to a regular map $S \rightarrow C$ around $p$ the restricted map 
$$f|_{\pi^{-1}(p)}: \pi^{-1}(p) \rightarrow C$$ 
is dominant (since $C$ is one-dimensional). 

Now the local monodromy around $p$ with respect to $S \dashrightarrow C$ is the same as the local monodromy around $\pi^{-1}(p)$ with respect to $S' \rightarrow C$. This local monodromy is the image of $\pi_1(V) \rightarrow \pi_1(C \setminus \cup_i p_i)$ where $V$ is some analytic open set whose analytic closure contains $\pi^{-1}(p)$. Now we can replace $V$ by a bigger open subset $V'$ which contains $f^{-1}(C \setminus \cup_i p_i) \cap \pi^{-1}(p)$ such that $V \subset V'$ is dense. Since the map $\pi_1(V) \rightarrow \pi_1(V')$ is surjective the images of 
$$\pi_1(V) \rightarrow \pi_1(C \setminus \cup_i p_i) \text{ and } \pi_1(V') \rightarrow \pi_1(C \setminus \cup_i p_i)$$ 
will be the same. Now if the local monodromy around $\pi^{-1}(p)$ were virtually abelian then the image of 
\begin{equation}\label{eq2}
\pi_1(V') \rightarrow \pi_1(C \setminus \cup_i p_i)
\end{equation}
would be virtually abelian. But this means that the image of  
$$\pi_1(f^{-1}(C \setminus \cup_i p_i) \cap \pi^{-1}(p)) \rightarrow \pi_1(C \setminus \cup_i p_i)$$
would be virtually abelian since it is a subgroup of the image in (\ref{eq2}). 

Now 
$$f^{-1}(C \setminus \cup_i p_i) \cap \pi^{-1}(p) \rightarrow C \setminus \cup_i p_i$$ 
is a dominant map. So by Lemma \ref{lem3} the induced image of fundamental groups is a finite index subgroup of $\pi_1(C \setminus \cup_i p_i)$ (which is not virtually abelian by the key observation above). Thus the local monodromy around $\pi^{-1}(p)$ cannot be virtually abelian. 

\begin{lem}\label{lem3} Let $\pi: S' \rightarrow S$ be a dominant morphism between irreducible varieties. If $S$ is normal then the image of the induced map $\pi_1(S') \rightarrow \pi_1(S)$ is a finite index subgroup. 
\end{lem}
\begin{proof}
See Lemma 3.3 (resp. Lemma 11) of \cite{k}. 

Alternatively, here is a short proof suggested by the referee. By taking general hyperplane sections of $S'$ we can find $S'' \hookrightarrow S'$ such that the composition $f: S'' \rightarrow S$ generically finite. Then we can restrict $f$ to some open $U \subset S$ such that $f: f^{-1}(U) \rightarrow U$ is finite \'etale. So the image $\pi_1(f^{-1}(U)) \rightarrow \pi_1(U)$ has finite index. Since $\pi_1(U) \rightarrow \pi_1(S)$ is surjective (this is where we use that $S$ is normal) this means the image of $\pi_1(f^{-1}(U)) \rightarrow \pi_1(S)$ has finite index. Since this map factors through $\pi_1(S') \rightarrow \pi_1(S)$ the result follows. 
\end{proof}

\subsubsection{The AME property for stacks}

The definition of the AME property can be extended to stacks as follows. We denote stacks using caligraphic font and their coarse moduli spaces by ordinary font. Our stacks will be integral, separated, normal Deligne-Mumford (DM) stacks of finite type over $\C$. Notice that by Keel-Mori \cite{km} a DM stack always has a coarse moduli scheme. 

We will need to use the (topological) fundamental group of a stack. For a reference see \cite{n}. If the stack is a smooth DM stack then this agrees with the orbifold fundamental group introduced by Thurston (see also \cite{alr}). 

If $\sX$ is a stack and $U \rightarrow \sX$ a regular morphism then we define the local monodromy around $p \in S \setminus U = D$ to be 
$$\mbox{Im} \left( \pi_1(V \cap U) \rightarrow \pi_1(\sX) \right)$$ 
where $V$ is a sufficiently small neighbourhood of $p$ as before. 

Then a variety $\oX$ containing the coarse space $X$ of $\sX$ as an open subset is an AME compactification of $\sX$ if the morphism $U \rightarrow \sX$ extends to a morphism $S \rightarrow \oX$ in a neighbourhood of $p$ whenever the local monodromy around $p$ is virtually abelian. As in Lemma \ref{lem1} this implies $\oX$ is complete.

\begin{Remark}
Note that we only consider a compactification $\oX$ of $X$ rather than $\sX$. There are basically two reasons to do this.
\begin{enumerate}
\item Under this definition the pair $(\sM_g, \oM_g)$ has the AME property (Theorem \ref{thm6}). On the other hand, it is not true that every morphism $C^* \rightarrow \sM_g$ from a punctured, smooth curve extends to a morphism $C \rightarrow \soM_g$ to the fine moduli stack (the central fibres might not be a stable curves). Since the local fundamental group of any puncture in $C^*$ is $\Z$ (and in particular abelian) this means that $(\sM_g, \soM_g)$ cannot have the AME property in any reasonable sense. 

Similarly, the pair $(\sA_g, A_g^{BB})$ has the AME property (Theorem \ref{thm5}) even though the Baily-Borel compactification $A_g^{BB}$ only compactifies $A_g$. So in both these cases it makes sense to consider compactifications of $X$ rather than $\sX$. 
\item Suppose $\soX$ is a compactification of $\sX$ with $\oX$ its coarse moduli scheme. Further suppose that we have a morphism $S \rightarrow \oX$ extending $U \rightarrow \sX \rightarrow X$. On $S$ the fibre product $\soX \times_{\oX} S \rightarrow S$ is (\'etale locally) of the form $[\tS/G] \rightarrow S$ where $\tS \rightarrow S$ is a finite cover and $G$ is the group of deck transformations. Thus $S \rightarrow \oX$ lifts to a morphism $\tS \rightarrow \soX$ on a finite cover $\tS$ of $S$. This means that the extra stack structure on $\soX$ can be ignored if you allow finite base changes and \'etale localization. 
\end{enumerate}

One might then be tempted to consider only schemes. However, the pair $(M_g, \oM_g)$ does {\em not} have the AME property because $M_g$ is simply connected (once you loose the stack structure on $\sM_g$ the fundamental group becomes trivial). So if we ignored stacks altogether we would lose some nice geometric examples of AME pairs like $(\sM_g, \oM_g)$ and $(\sA_g, A_g^{BB})$. 
\end{Remark}

\section{Properties of AME compactifications}\label{s2}

\begin{prop}\label{prop4} Let $\sX$ and $\sY$ be two DM stacks whose coarse spaces $X$ and $Y$ have normal compactifications $\oX$ and $\oY$. Let $\pi: \sY \rightarrow \sX$ be a morphism which extends to a finite morphism $\pi: \oY \rightarrow \oX$. If $(\sX,\oX)$ has the AME property then $(\sY,\oY)$ has the AME property. Conversely, if $\pi: \sY \rightarrow \sX$ is also representable, finite, \'etale and $(\sY,\oY)$ has the AME property then $(\sX,\oX)$ has the AME property. 
\end{prop}
\begin{proof}
Suppose $(\sX,\oX)$ has the AME property and consider a morphism $U \rightarrow \sY$ from an open subset $U \subset S$. If the local monodromy around $p \in S \setminus U$ with respect to $U \rightarrow \sY$ is virtually abelian then the local monodromy around $p$ with respect to $U \rightarrow \sX$ is the image of the composition $\pi_1(U) \rightarrow \pi_1(\sY) \rightarrow \pi_1(\sX)$ and is also virtually abelian. Hence the composition morphism $U \rightarrow \sX$ extends to a regular morphism $S \rightarrow \oX$. Since $\pi: \oY \rightarrow \oX$ is finite the morphism $S \rightarrow \oX$ lifts to a morphism $S \rightarrow \oY$ by Lemma \ref{lem5}. This shows $(\sY,\oY)$ has the AME property.

Conversely, suppose that $(\sY,\oY)$ has the AME property and $\pi: \sY \rightarrow \sX$ is \'etale. Consider a morphism $U \rightarrow \sX$ with virtually abelian local monodromy around $p \in S \setminus U$. Consider the fibre product 
\begin{equation*}
\begin{CD}
\tS @<{\supset}<< \widetilde{U} = U \times_{\sX} \sY @>>> \sY \\
@VVV @VVV @V{\pi}VV \\
S @<{\supset}<< U @>>> \sX 
\end{CD}
\end{equation*}
and denote by $\tS$ the normal closure of $S$ in the function field of $\widetilde{U}$. For any point $\widetilde{p}$ in the preimage of $p$ under $\tS \rightarrow S$ the local monodromy around $\widetilde{p}$ with respect to $\widetilde{U} \rightarrow \sX$ is virtually abelian. Since $\pi: \sY \rightarrow \sX$ is finite, \'etale the map $\pi_1(\sY) \rightarrow \pi_1(\sX)$ is injective so the local monodromy around $\widetilde{p}$ with respect to $\widetilde{U} \rightarrow \sY$ must also be virtually abelian. This means the morphism $\widetilde{U} \rightarrow \sY$ extends to a morphism $\tS \rightarrow \oY$ in a neighbourhood of $\widetilde{p}$. Since this is true for every such $\widetilde{p}$ the morphism $\widetilde{U} \rightarrow \sX$ extends to a morphism $\tS \rightarrow \oX$ in a neighbourhood of the preimage of $p$ under $\tS \rightarrow S$. 

If, in the notation of Lemma \ref{lem4}, we look at $\tS \rightarrow S \dashrightarrow \oX$ then the composition is regular so $S \rightarrow \oX$ is regular. Thus $U \rightarrow \sX$ also extends to a morphism $S \rightarrow \oX$ in a neighbourhood of $p$. This shows $(\sX, \oX)$ has the AME property. 

\begin{lem}\label{lem5} Consider the composition 
\[
\begin{xy} 
\xymatrix{ 
S \ar@{-->}[r]^h & Y \ar[d]^{\pi} \\
 & \sX } 
\end{xy} 
\]
where $S$ and $Y$ are varieties, $S$ is normal and $\sX$ is a DM stack. If $\pi$ is a finite morphism then the rational morphism $h$ is regular if and only if $\pi \circ h$ is regular.
\end{lem}
\begin{proof}

As the referee points out, this actually follows from Lemma \ref{lem4}. We need to show that if $\pi \circ h$ is regular then so is $h$. If one takes $\tS := S \times_\sX Y$ then we are left with the diagram 
\[
\begin{xy}
\xymatrix{
\tS \ar[d]_{\pi'} \ar[dr] & \\
S \ar@{-->}[r]^{h} & Y
}
\end{xy}
\]
where $\pi'$ is a finite morphism. So $h$ is regular by Lemma \ref{lem4}. 
\end{proof}
\end{proof}

\begin{cor}\label{cor8} Suppose $(\sX,\oX)$ has the AME property and let $\sY \subset \sX$ be a closed, normal DM substack. Denote by $X$ and $Y$ the correponding coarse moduli schemes. Then $Y \subset X$ and if we denote by $\oY$ the normalization of the closure of $Y$ in $\oX$ then the pair $(\sY,\oY)$ also has the AME property. 
\end{cor}
\begin{proof}
This follows from Proposition \ref{prop4} by considering a closed immersion $\sY \hookrightarrow \sX$. We just need to show that $\sY \subset \sX$ implies $Y \subset X$. 

To do this we use that a DM stack $\sX$ over $\C$ (or more generally a tame stack in the sense of \cite{aov}) is \'etale locally $\sX$ of the form $[T/G]$ for some scheme $T$ with an action of a finite group $G$ (see, for example, the introduction of \cite{aov}). More precisely, this means that there exists a finite \'etale cover $\tX \rightarrow X$ of the coarse moduli scheme of $\sX$ such that $\sX \times_X \tX \rightarrow \tX$ is locally in the Zariski topology isomorphic to $[T/G] \rightarrow (T/G)$. 

The fact that $\sY \hookrightarrow \sX$ is some closed substack means that 
$$\sY \times_X \tX = [T'/G] \hookrightarrow [T/G] = \sX \times_X \tX$$
where $T' \hookrightarrow T$ is a closed $G$-invariant subscheme. Now consider the following commutative diagram 
\[
\begin{xy}
\xymatrix{
\sY \ar[d] & \sY \times_X \tX = [T'/G] \ar[d] \ar[r] \ar[l] & (T'/G) = Y \times_X \tX \ar[r] \ar[d] & Y \ar[d] \\
\sX \ar[d] & \sX \times_X \tX = [T/G] \ar[d] \ar[r] \ar[l] & (T/G) = \tX \ar[r] & X \\
X & \tX = (T/G) \ar[l] & & 
}
\end{xy}
\]
Looking at the right-most square we see that $T'/G \hookrightarrow T/G$ because $T' \hookrightarrow T$. Since $\tX \rightarrow X$ is an \'etale cover it follows $Y \hookrightarrow X$. 
\end{proof}

\begin{prop}\label{prop6} If $(\sX_1, \oX_1)$ and $(\sX_2, \oX_2)$ are two AME pairs then $(\sX_1 \times \sX_2, \oX_1 \times \oX_2)$ is also an AME pair.
\end{prop}
\begin{proof}
Consider an open subvariety $U \subset S$ of a normal variety $S$ together with a morphism $U \rightarrow \sX_1 \times \sX_2$ with virtually abelian monodromy around $p \in S \setminus U$. The composition with the two projection maps yields morphisms $U \rightarrow \sX_1$ and $U \rightarrow \sX_2$. The monodromy around $p$ with respect to these maps is also virtually abelian. Since $\oX_1$ and $\oX_2$ are AME compactifications these maps extend to regular morphisms $S \rightarrow \oX_1$ and $S \rightarrow \oX_2$ in a neighbourhood of $p$. This gives a regular morphism $S \rightarrow \oX_1 \times \oX_2$. 

\begin{Remark} Note that the coarse scheme of a product of two DM stacks is the product of the coarse schemes. To see this suppose $\sX_1$ and $\sX_2$ are DM stacks with coarse schemes $X_1$ and $X_2$ and let $Z$ denote the coarse scheme of $\sX_1 \times \sX_2$. We must show that $Z = X_1 \times X_2$. 

Now \'etale locally $\sX_1$ and $\sX_2$ are both quotients of a scheme by a finite group \cite{aov}. This means that Zariski locally on $X_1$ and $X_2$ we can find \'etale covers $\tX_1 \rightarrow X_1$ and $\tX_2 \rightarrow X_2$ such that $\tX_1 \times_{X_1} \sX_1 \rightarrow \tX_1$ is isomorphic to $[T/G] \rightarrow (T/G)$ and $\tX_2 \times_{X_2} \sX_2$ is isomorphic to $[T'/G'] \rightarrow (T'/G')$ for some schemes $T,T'$ with finite group actions $G,G'$. 

Consider the following commutative diagram where every square is a fibre product:
\[
\begin{xy}
\xymatrix{
\sX_1 \times \sX_2 \ar[d] & W := (\sX_1 \times \sX_2) \times_{X_1 \times X_2} (\tX_1 \times \tX_2) \ar[d] \ar[l] \\
Z \ar[d]^f & Z \times_{X_1 \times X_2} (\tX_1 \times \tX_2) \ar[d]^g \ar[l] \\
X_1 \times X_2 & \tX_1 \times \tX_2. \ar[l]}
\end{xy}
\]
Now $\tX_1 \times \tX_2 = (T/G) \times (T'/G')$ is the coarse scheme of 
$$W = (\sX_1 \times_{X_1} \tX_1) \times (\sX_1 \times_{X_2} \tX_2) = [T/G] \times [T'/G'].$$
On the other hand, since $Z$ is the coarse scheme of $\sX_1 \times \sX_2$ it must be that $Z \times_{X_1 \times X_2} (\tX_1 \times \tX_2)$ is also the coarse scheme of $W$ (Corollary 3.3 of \cite{aov}). Thus $g$ must be an isomorphism and since $\tX_1 \times \tX_2 \rightarrow X_1 \times X_2$ is an \'etale cover $f$ must also be an isomorphism. 
\end{Remark}
\end{proof}

\Example We use \ref{cor8} and \ref{prop6} to construct more examples of AME compactifications. Choose four general points $p_1, \dots, p_4$ on $\P^2$ and let $l_{ij}$ denote the line between $p_i$ and $p_j$. If we take $W_2 = \P^2 \setminus \cup_{1 \le i < j \le 4} l_{ij}$ and $\oW_2 = \P^2$ then $(W_2,\oW_2)$ is an AME pair. To see this consider the morphisms $\pi_i: W_2 \rightarrow \P^1 \setminus \{0,1,\infty\}$ for $1 \le i \le 4$ which take a point $q \in W_2$ and map it to the line $[qp_i]$. It is not hard to see that the product of these four maps gives an embedding
$$\pi_1 \times \dots \times \pi_4: W_2 \rightarrow (\P^1 \setminus \{0,1,\infty\})^{\times 4} \subset (\P^1)^{\times 4}$$
whose closure we denote $\oW_2$. Actually, $\oW_2$ is nothing but $\P^2$ blown up at the four points $p_i$. Since $(\P^1 \setminus \{0,1,\infty\},\P^1)$ is an AME pair so is $(W_2,\oW_2)$ by an application of \ref{cor8} and \ref{prop6}. In fact $(W_2,\oW_2) \cong (M_{0,5},\oM_{0,5})$. Since $\oW_2$ projects onto $\P^2$ we find that $(W_2,\P^2)$ is also an AME pair. 

More generally, we may take $n+2$ general points $p_1, \dots, p_{n+2} \in \P^n$ and remove all hyperplanes through any $n$ of them to obtain $W_n$. As before, we have morphisms $\pi_i: W_n \rightarrow W_{n-1}$ for $1 \le i \le n+2$ taking a point $q \in X$ to the line $[qp_i]$. The product of these maps 
$$\iota_n = \prod_{i=1}^{n+2} \pi_i: W_n \rightarrow (W_{n-1})^{\times n+2} \subset (\P^{n-1})^{n+2}$$
gives an embedding of $W_n$. Denote by $\oW_n$ the closure of $\iota_n(W_n)$ inside $(\P^{n-1})^{n+2}$. By induction, $(W_{n-1}, \P^{n-1})$ is an AME pair and so $(W_n, \oW_n)$ is also an AME pair. On the other hand, $\iota_n$ extends to an embedding of $\P^n \setminus \cup_i p_i$ while it is undefined at the points $\cup_i p_i$. If one blows up $\P^n$ at the points $\cup_i p_i$ then $\iota_n$ extends and we find that $\oW_n$ is precisely this blowup. Thus $\oW_n$ projects onto $\P^n$ and so $(W_n, \P^n)$ is an AME pair. 

\begin{cor}\label{cor5} Let $\oX_1$ and $\oX_2$ be two AME compactifications of $\sX$. Then there exists an AME compactification $\oX$ of $\sX$ which dominates $\oX_1$ and $\oX_2$.
\end{cor}
\begin{proof}
$\oX$ dominates $\oX_1$ and $\oX_2$ if it admits regular maps $\pi_1: \oX \rightarrow \oX_1$ and $\pi_2: \oX \rightarrow \oX_2$ which extend the identity map on $X$. Take $\oX$ to be the normalization of the closure of $X$ embedded inside $\oX_1 \times \oX_2$ by the diagonal map. The two projections give us $\pi_1$ and $\pi_2$ so it remains to show $(\sX, \oX)$ has the AME property. By Proposition \ref{prop6} $(\sX_1 \times \sX_2, \oX_1 \times \oX_2)$ is an AME pair so by Corollary \ref{cor8} the closure $\oX$ of $X$ in $\oX_1 \times \oX_2$ is an AME compactification.
\end{proof}

\begin{cor}\label{cor6} If $\sX$ has the AME property then it has a unique, maximal AME compactification $\oX_{ame}$ in the sense that for any other AME compactification $\oX$ there is a birational, regular map $\oX_{ame} \rightarrow \oX$. 
\end{cor}
\begin{proof}
By Corollary (\ref{cor5}) any two AME compactifications of $\sX$ are dominated by a larger AME compactification. If there is no maximal AME compactification there exists a sequence $\oX_0 \xleftarrow{\pi_0} \oX_1 \xleftarrow{\pi_1} \oX_2 \xleftarrow{\pi_2} \dots$ of AME compactifications where each birational map $\pi_i: \oX_{i+1} \rightarrow \oX_i$ between complete, normal, reduced schemes is not an isomorphism. Take any normal crossing compactification $\oY$ of $X$. The local fundamental group of a neighbourhood of $p \in \oY \setminus X$ is free abelian. This means every open inclusion $X \hookrightarrow \oX_i$ extends to a regular map $f_i: \oY \rightarrow \oX_i$. By Lemma (\ref{lem6}) we know $\oX_n \cong \oX_{n+1}$ for $n \gg 0$ so this is impossible. 

If $\oX$ and $\oX'$ are two maximal AME compactifications then there are birational, regular maps $\oX \rightarrow \oX'$ and $\oX' \rightarrow \oX$ which implies $\oX \cong \oX'$ and so a maximal AME compactification is unique. 

\begin{lem}\label{lem6} Let 
$X_0 \xleftarrow{\pi_0} X_1 \xleftarrow{\pi_1} X_2 \xleftarrow{\pi_2} \dots$
be a sequence of birational morphisms between complete, normal, reduced schemes. Let $Y$ be a complete, reduced scheme equipped with dominant morphisms $f_i: Y \rightarrow X_i$ satisfying $\pi_i \circ f_{i+1} = f_i$ for each $i$. Then $X_n \cong X_{n+1}$ for $n \gg 0$.
\end{lem}
\begin{proof}
Here is a short proof suggested by the referee. We have natural inclusions $Y \times_{X_{n+1}} Y \hookrightarrow Y \times_{X_n} Y$. These form a decreasing sequence of closed subsets $Y \times_{X_n} Y$ inside the Noetherian scheme $Y \times_{X_0} Y$. Thus they must stabilize at some point. 

Now if $\pi_n: X_{n+1} \rightarrow X_n$ is not an isomorphism then by Zariski's main theorem there must be two points $q_1 \ne q_2 \in X_{n+1}$ such that $\pi_n(q_1) = \pi_n(q_2)$. Since the $f_i$ are surjective (they are dominant and all varieties are complete) we can find points $p_1 \ne p_2$ such that $f_n(p_i) = q_i$. But then $(p_1,p_2)$ belongs to $Y \times_{X_n} Y$ but not to $Y \times_{X_{n+1}} Y$. Hence $Y \times_{X_n} Y \ne Y \times_{X_{n+1}} Y$. 

Since the subschemes $Y \times_{X_n} Y \subset Y \times_{X_0} Y$ are the same for $n \gg 0$ this means that $\pi_n$ is an isomorphism for $n \gg 0$. 
\end{proof}
\end{proof}

\begin{Remark}
By construction, any AME compactification $\oX$ of $\sX$ is obtained from $\oX_{ame}$ by contracting parts of the boundary. Conversely, given a contraction of the boundary $\oX_{ame} \rightarrow \oX$, $(\sX, \oX)$ is an AME pair. So to understand all AME compactifications of $\sX$ it is enough to understand the maximal one $\oX_{ame}$ and all possible contractions of its boundary. 
\end{Remark}

\begin{prop}\label{prop7} If $\sX$ has the AME property and $i: \sY \rightarrow \sX$ is a locally closed embedding then there exists a regular morphism $\oY_{ame} \rightarrow \oX_{ame}$ which extends $i: Y \rightarrow X$.
\end{prop}
\begin{proof}
The closure $\oY$ of $Y$ in some $\oX_{ame}$ is an AME compactification. By the defining property of $\oY_{ame}$ there exist morphisms $\oY_{ame} \rightarrow \oY \rightarrow \oX_{ame}$ whose composition extends $i: Y \rightarrow X$
\end{proof}

The following theorem shows that varieties having the AME property are numerous. 

\begin{prop} Let $X \subset \oX$ be a dense, open immersion with $\oX$ normal, complete variety. Then there exists an open $X^o \subset X$ such that $(X^o, \oX)$ has the AME property. 
\end{prop}
\begin{proof}
Suppose $\oX \subset \P^n$ is projective. By example 2 there exists an open subvariety $H_n \subset \P^n$ such that $(H_n,\P^n)$ is an AME pair. Let $X^o = X \cap H_n$. The closure of $X^o$ is $\oX$ so that $(X^o,\oX)$ is an AME pair by Corollary \ref{cor8}. 

If $\oX$ is not projective then by Chow's lemma there exists a proper, birational morphism $\pi: \oY \rightarrow \oX$ where $\oY$ is normal and projective. Then by the argument above we can find $Y^o \subset \oY$ such that $(Y^o, \oY)$ is an AME pair. Restricting $Y^o$ further we can even assume $Y^o \rightarrow \pi(Y^o)$ is an isomorphism with $\pi(Y^o) \subset X$. This immediately implies $(X^o, \oX)$ is an AME pair where $X^o = \pi(Y^o)$.
\end{proof}

We end the section with a Lemma, which though not entirely necessary, will simplify some of the subsequent proofs.

\begin{lem}\label{lem7} To show $(\sX,\oX)$ is an AME pair it suffices to check the extension property for normal surfaces $U \subset S$. 
\end{lem}
\begin{proof}
Suppose $U \rightarrow \sX$ does not extend to a regular map $S \rightarrow \oX$ in a neighbourhood of $p \in S \setminus U$. Assuming $\dim(S) > 2$, we construct a codimension one subvariety $T \subset S$ through $p$ with the following property. Denote by $\pi: \tT \rightarrow T$ the normalization of $T$. Then $\pi^{-1}(U \cap T) \rightarrow \sX$ does not extend to a regular map $\tT \rightarrow \oX$ in a neighbourhood of some point $\tp \in \pi^{-1}(p)$. If the local monodromy around $p$ with respect to $U \rightarrow \sX$ is some virtually abelian group $G$ then the local monodromy with respect to $\pi^{-1}(U \cap T) \rightarrow \sX$ around $\tp$ is a subgroup of $G$ and hence also virtually abelian. Iterating we arrive at a surface which violates the AME property of $(\sX, \oX)$, thus proving the lemma. 

Consider the closure $S'$ of the graph of $S \dashrightarrow \oX$ in $S \times \oX$. Since $\oX$ is complete the projection $\pi_1: S' \rightarrow S$ is a proper, birational morphism. If the fibre $\pi_1^{-1}(p)$ is zero dimensional then, by Zariski's main theorem (EGA III Corollary 4.4.9), in a neighbourhood of $p$, $S \cong S'$ implying $S \rightarrow \oX$ is regular (contradiction). Let $C \subset \pi_1^{-1}(p)$ be a curve. Locally around some point in $C$ take an irreducible, codimension one normal subscheme containing a generic point of $C$ and passing through a point of $S'$ not on the exceptional locus of $\pi_1: S' \rightarrow S$. This is possible if $\dim(S') > 2$. Denote its closure in $S'$ by $T' \subset S'$ and let $T = \pi_1(T) \subset S$. 

Then $T$ is an irreducible, codimension one subscheme of $S$ which passes through $p$ and intersects $U$. Moreover, $T \cap U \rightarrow \sX$ does not extend over $p$ since the closure of the graph of $T \dashrightarrow \oX$ inside $S \times \oX$ contains $C$ in the fibre over $p \in T$. In fact, for the same reason $\tT \dashrightarrow \oX$ is not regular in a neighbourhood of $\pi^{-1}(p)$ and hence in a neighbourhood of some $\tp \in \pi^{-1}(p)$ as desired. 
\end{proof}

\begin{Remark} Instead of working with arbitrary varieties one can restrict to the category of projective varieties and every result in this section still holds. Instead of a unique maximal AME compactification $\oX_{ame}$ of $\sX$ there is a unique maximal projective AME compactification $\oX_{amep}$ and a proper morphism $\oX_{ame} \rightarrow \oX_{amep}$. 
\begin{prob} Find an example where $\oX_{ame} \not\cong \oX_{amep}$. 
\end{prob}
\end{Remark}

\section{The Moduli Space of Curves}\label{s3}

Denote by $\sM_{g,n}$ the DM moduli stack of stable, smooth, $n$-pointed, genus $g$ curves and by $\soM_{g,n}$ its Deligne-Mumford compactification (as a DM stack). Denote by $M_{g,n}$ and $\oM_{g,n}$ the corresponding coarse spaces. See \cite{hm} and \cite{v} for an introduction to some of the theory involving these spaces.

In this section we prove:

\begin{thm}\label{thm6}
$\oM_{g,n}$ is the maximal AME compactification of $\sM_{g,n}$. 
\end{thm}

\subsection{Some Results on Monodromy in Families of Stable Curves}\label{s31}

\subsubsection{Definitions: Local and Global Monodromy of a Stable Family}

As before, let $S$ be a normal variety and $U \subset S$ an dense, open subset. By a family of stable curves $\pi: \Cu \rightarrow U$ we mean a flat family of stable $n$-pointed, genus $g$ curve. We denote the marked points by $\{p_i\}$ where $1 \le i \le n$.

If the general fibre of $\pi: \Cu \rightarrow U$ is smooth then we get a morphism $U' \rightarrow \sM_{g,n}$ from some dense, open $U' \subset U$. Let $T \subset S$ be a connected, reduced, complete subscheme. As in section \ref{s2.1} we restrict $U'$ to a smaller open subset which is disjoint from $T$ and any singularities of $S$. After blowing up we can assume the complement of $U'$ is a normal crossing divisor. Then we define the {\em local monodromy around $T$} as the image 
$$\mbox{Im}(\pi_1(U' \cap V) \rightarrow \pi_1(\sM_{g,n}))$$
where $V$ is a small neighbourhood of $T$ as in section \ref{s2.1}. Similarly, if $S$ is complete, the {\em global monodromy} is the image of $\pi_1(U') \rightarrow \pi_1(\sM_{g,n})$. 

This is precisely the definition from section \ref{s2.1} applied to the morphism $U' \rightarrow \sX$ when $\sX = \sM_{g,n}$. If the general fibre of $\pi: \Cu \rightarrow U$ is not smooth then we can find a dense, open $U' \subset U$ over which $\Cu \setminus \{p_i\}$ is of constant topological type. This gives a map $U' \rightarrow \sM' \subset \soM_{g,n}$ into some boundary stratum. We then define the {\em local monodromy} to be the image
\begin{equation}\label{eq1}
\mbox{Im}(\pi_1(U' \cap V) \rightarrow \pi_1(\sM'))
\end{equation}
as above. Similarly, the {\em global monodromy} is the image of $\pi_1(U') \rightarrow \pi_1(\sM')$. 

Let us study a little more carefully this map $U' \rightarrow \sM'$. We restrict our attention to the family $\pi': \Cu' := \Cu|_{U'} \rightarrow U'$. After possibly having to pull back $\Cu'$ to a finite cover of $U'$, the normalization of $\Cu'$ is the disjoint union $\Cu_1 \sqcup \dots \sqcup \Cu_m$ where each $\Cu_i \rightarrow U'$ is a family of smooth, marked, stable curves. By Lemma \ref{lem0} below $\Sing(\pi') \rightarrow S$ is finite, unramified and surjective. So, after a finite \'etale base change $\tU' \rightarrow U'$, the singular points becomes sections (see for instance EGA IV, (18.4.7)). Hence $\Cu'$ is obtained from $\Cu_1 \sqcup \dots \sqcup \Cu_m$ by glueing pairs of these ``special'' sections (``special'' in order to distinguish them from the sections induced by the marked points). Thus we get a map 
$$\tU' \rightarrow \sM_{g_1,n_1+n_1'} \times \dots \times \sM_{g_m,n_m+n_m'}$$
where $g_i$ is the genus of the fibres of $\Cu_i$, $n_i$ is the number of marked points of $\Cu'$ lying on $\Cu_i$ and $n_i'$ is the number of special sections on $\Cu_i$. Note that $n = \sum_i n_i$ while $g = \sum_i g_i + \frac{1}{2} \sum_i n_i' - (m-1)$. 

Now we could try to define the local monodromy around $T \subset S$ as
$$\mbox{Im} \left( \pi_1(\tV \cap \tU') \rightarrow \pi_1(\sM_{g_1,n_1+n_1'} \times \dots \times \sM_{g_m,n_m+n_m'}) \right)$$ 
where $\tV$ is a small neighbourhood of the preimage of $T$ in $\tS$ (where $\tS$ is the normal closure of $S$ in the function field of $\tU'$). There was a choice of cover $\tU'$ so this image is only defined up to the equivalence $\sim$ where $G_1 \sim G_2$ if they are contained in a common group $G$ as subgroups of finite index. Fortunately, this is sufficient for us since we get a well defined concept of the image being virtually abelian. Even better, modulo $\sim$ this image is equal to the local monodromy defined by (\ref{eq1}). 

\begin{lem}\label{lem0} Let $\pi': \Cu' \rightarrow U'$ be a family of stable curves over (an irreducible) variety $U'$. Then $\Sing(\pi') \subset \Cu'$ has a natural closed subscheme structure which commutes with base change and is finite unramified over $U'$. 
\end{lem}
\begin{proof}
See \cite{jo} Lemma 4.3 or \cite{mb} Lemma 3 or \cite{dm}. 
\end{proof}

\subsubsection{Virtually Abelian Monodromy}

In what follows we will often identify the mapping class group $\Gamma_{g,n} = \pi_1(\sM_{g,n})$ with the outer automorphism group of the fundamental group of a genus $g$ Riemann surface with $n$ punctures. It is well known that $\Gamma_{g,n}$ is generated by Dehn twists. More generally, we also consider $\Gamma_{g,n}(m) := \mbox{ker}(\Gamma_{g,n} \rightarrow Sp_{2g}(\Z/m\Z))$ which is the mapping class group of curves with a level $m$ structure. 

See \cite{d} for a short introduction to vanishing cycles and Dehn twists. 

\begin{prop}\label{prop3} Let $S$ be a normal variety and suppose $\pi: \Cu \rightarrow S$ is a (semi)stable family of pointed curves. Let $T \subset S$ be a connected, integral, complete subscheme. If the local monodromy around $T$ is virtually abelian then the global monodromy of $\Cu|_T \rightarrow T$ is also virtually abelian. 
\end{prop}
\begin{proof}
Since we only deal with virtually abelian monodromy it suffices to prove the result for the pullback of $\Cu$ to any finite branched cover of $S$. From the definition of monodromy above we can assume the general fibre of $\Cu \rightarrow S$ is smooth (even though the general fibre over $T$ may still be singular) and that the monodromy around $T$ is abelian and lies in $\Gamma_{g,n}(m)$ for some $m \ge 3$.

The local monodromy around $T$ is the image of $\pi_1(U \cap V) \rightarrow \pi_1(\sM_{g,n})$ where $V$ is a small neighbourhood of $T$. If the general fibre over $T$ is also smooth then we can extend $U$ to include an open dense subset of $T$. But then the global monodromy on $T$ is a subgroup of the local monodromy around $T$ and thus must be abelian (and we are done).

If the general fibre over $T$ is singular consider an arc $A$ joining a general point of $S$ to a general point $0 \in T$. The restriction of $\Cu$ to $A$ is a family of smooth, pointed curves degenerating to a semistable curve over $0 \in A$. Consider the vanishing cycles $v_1, \dots, v_k$ of this one dimensional family and choose simple disjoint loops (which we also call $v_1, \dots, v_k$) representing them. Thus these $v_i$ shrink to a point as we approach $0 \in A$. Then if we cut $\Cu|_A$ along these loops we are left with a family of (possibly disconnected) Riemann surfaces with boundary except now all the fibres have the same topological type while the fibre over $0 \in A$ is a punctured Riemann surface. We can similarly do this in a small analytic neighbourhood of $A$. Now we would like to do this globally. 

The local monodromy around $T$ is abelian, lies in $\Gamma_{g,n}(m)$ and contains the Dehn twist in the multi-loop $\sum_i a_i v_i$ (for some $a_i \in \Z^{>0}$). By Theorem \ref{thm2} this means that the local monodromy around $T$ preserves the vanishing cycles (we do not actually need the full strength of Theorem \ref{thm2}). This is what allows us, in a neighbourhood of $T$, to consistently choose simple, disjoint loops representing $v_1, \dots, v_k$. 

If we now again cut the fibres of $\Cu$ along these loops we are left with a family of (possibly disconnected) Riemann surfaces with boundary but now the general fibre over $S$ is of the same topological type to a general fibre over $T$ which is a punctured Riemann surface. Since the monodromy around $T$ is abelian the monodromy over $T$ must also be abelian (using the same argument as above). 
\end{proof}

\begin{prop}\label{prop8} A family of pointed, stable curves $\Cu \rightarrow S$ over a normal variety $S$ is isotrivial if the global monodromy is virtually abelian. 
\end{prop}
\begin{proof}
Since any two points in the base can be connected by a series of irreducible curves it suffices to prove the result when $S$ is a curve.  If $S$ is not smooth or complete take its normalization and compactify. Taking a branched cover of this new curve the pullback of $\Cu$ will extend to a family of stable curves so that we may assume from now on that $S$ is a smooth, complete curve. 

By Lemma \ref{lem0} it suffices to prove the result when the general fibre of $\Cu$ is smooth. We first deal with the case when $\Cu \rightarrow S$ is unpointed. Denote by $U \subset S$ the largest open subset over which the fibres of $\Cu$ are smooth. The first cohomology group $H^1$ of the fibres of $\Cu|_U \rightarrow U$ gives a variation of principally polarized Hodge structures of weight one. The induced monodromy map on this variation is the composition $\pi_1(U) \rightarrow \Gamma_g \rightarrow Sp_{2g}(\Z)$ so that the monodromy of this variation is virtually abelian. By Corollary 4.2.9 of \cite{del}, the connected component containing the identity of the Zariski closure of the image of the map $\pi_1(U) \rightarrow Sp_{2g}(\C)$ is semi-simple. Since the image of $\pi_1(U) \rightarrow Sp_{2g}(\Z)$ is (virtually) abelian so is its Zariski closure in $Sp_{2g}(\C)$. Since a connected semisimple abelian group is trivial the monodromy must be finite and hence the variation of Hodge structure is isotrivial. Consequently, by Torelli's theorem, $\Cu|_U \rightarrow U$ is isotrivial. 

To deal with the pointed case consider the short exact sequence of groups 
$$1 \rightarrow \pi_1(C - \{n \mbox{ points}\}) \rightarrow \Gamma_{g,n+1} \rightarrow \Gamma_{g,n} \rightarrow 1$$
obtained by forgetting the $(n+1)$st point of the curve $C$ (if $C$ is a curve of genus $g$ then $n \ge 3$ if $g=0$ and $n \ge 1$ if $g=1$). We complete the proof by induction on $n$. The argument above proves the base case. Note that the base case $(g,n)=(1,1)$ is also covered by the argument above while the base case $(g,n)=(0,3)$ is obvious since any such family is trivial. 

Suppose $\Cu \rightarrow S$ is a stable family with $n+1$ marked points given by sections $\sigma_1, \dots, \sigma_{n+1}$. Since the global monodromy $G \subset \Gamma_{g,n+1}$ of $\Cu$ is abelian the image of $G$ in $\Gamma_{g,n}$ is also abelian so by induction (and after a finite base change if necessary) $\Cu$ is a trivial family $C \times S$ with $n$ constant sections corresponding to points $p_1, \dots, p_n \in C$. The $(n+1)$st section is given by a morphism $f: S \rightarrow C$. 

The kernel of $G \rightarrow \Gamma_{g,n}$ is the image of $f_\ast: \pi_1(S \setminus f^{-1}(\cup_i p_i)) \rightarrow \pi_1(C \setminus \cup_i p_i)$ which is therefore also abelian. By Lemma \ref{lem3}, if $f$ is dominant the image of $f_\ast$ has finite index in $\pi_1(C \setminus \cup_i p_i)$. Since there are no abelian subgroups of finite index inside $\pi_1(C \setminus \cup_i p_i)$ the map $f$ must be constant. This completes the proof.
\end{proof}

\begin{Remark}
It is possible to give an alternate proof of Proposition \ref{prop8} by using Theorem \ref{thm2} below in place of Corollary 4.2.9 from \cite{del}. The idea is that there is a nice classification of abelian subgroups of $\Gamma_{g,n}$ given by Ivanov who builds on the paper \cite{blm}. To explain this let $M$ be a surface of genus $g$ with $n$ punctures. A {\em system of circles} $C$ on $M$ is a collection of pairwise non-isotopic, non-intersecting circles on $M$. Denote by $M_C$ the surface obtained by cutting $M$ along $C$. A subgroup of $\Gamma_{g,n}(m)$ is a {\em $C$-subgroup} if it is generated by Dehn twists around (some of the) circles in $C$ and by (at most) one pseudo-Anosov map on each component of $M_C$ each of which fixes the boundary pointwise. The following theorem appears in \cite{iv} (page 205):

\begin{thm}\label{thm2} A subgroup of $\Gamma_{g,n}(m)$ ($m \ge 3$) is abelian if and only if it is a $C$-subgroup for some system of circles $C$.
\end{thm}
\end{Remark}

\subsection{Proof of Theorem \ref{thm6}: $\sM_{g,n}^{ame} = \oM_{g,n}$}\label{ss4}

As usual we have a normal variety $S$ and a dense, open $U \subset S$ with a morphism $U \rightarrow \sM_{g,n}$. Suppose the local monodromy around $p \in S$ is virtually abelian. We can assume $S$ is a small neighbourhood of $p$. We want to show that we get a regular map $S \rightarrow \oM_{g,n}$. 

Now $\soM_{g,n}$ has a finite cover $Z \rightarrow \soM_{g,n}$ by a projective scheme $Z$ (see \cite{kv} or \cite{l}). If we let $\tS$ be the normalization of the closure of $S \times_{\soM_{g,n}} Z$ inside $S \times Z$ then we get a projective, generically finite morphism $\tS \rightarrow S$ so that the composition $\tS \rightarrow \soM_{g,n}$ is a regular map resolving the rational map $S \dashrightarrow \soM_{g,n}$. This gives us a family of stable curves over $\tS$. 

Let $T$ be an irreducible component of the fibre $\tS_p$ over $p \in S$. The local monodromy around $T$ is virtually abelian so by Proposition \ref{prop3} the global monodromy on $T$ is also virtually abelian. Hence, by Proposition \ref{prop8}, the family over $T$ must be isotrivial. 

Let $\tS \xrightarrow{f} S'' \xrightarrow{g} S$ be the Stein factorization of $\tS \rightarrow S$ where $f: \tS \rightarrow S''$ has connected fibres and $g: S'' \rightarrow S$ is finite. Let $q \in g^{-1}(p)$ and consider the fibre $\tS_q$. By the above, the image of each irreducible component of $\tS_q$ is a point in $M_{g,n}$ and since $\tS_q$ is connected the whole fibre is collapsed to a point. This means $S'' \rightarrow M_{g,n}$ is regular in a neighbourhood of $g^{-1}(p)$ ($\tS$ and subsequently $S''$ are normal). Since $S$ is normal, Lemma \ref{lem4} implies $S \rightarrow M_{g,n}$ is also regular in a neighbourhood of $p$. This proves that $\oM_{g,n}$ is an AME compactification. 

To see that $\oM_{g,n}$ is maximal fix any $p \in \oM_{g,n}$ and consider the rational map $h: \oM_{g,n} \dashrightarrow \sM_{g,n}^{ame}$ in a neighbourhood of $p$. The map $\soM_{g,n} \rightarrow \oM_{g,n}$ \'etale locally around $p$ is of the form $[Y/G] \rightarrow Y/G$ where $G$ is a finite group and the boundary of $Y$ (i.e. the preimage of $\partial \oM_{g,n}$) is normal crossing (this is by definition what it means for $\soM_{g,n}$ to have normal crossing boundary). This means that the local fundamental group around any point on $Y$ is abelian so that we have a regular map $Y \rightarrow \sM_{g,n}^{ame}$. Thus we end up with the diagram
\[
\begin{xy}
\xymatrix{
Y \ar[d] \ar[drr] & & \\
\tV \ar[r] & V \ar@{-->}[r]^h & \sM_{g,n}^{ame}
}
\end{xy}
\]
where $V \subset \oM_{g,n}$ is a neighbourhood of $p$, $\tV = Y/G \rightarrow V$ is a finite \'etale cover and $Y \rightarrow \tV = Y/G$ is finite. It follows by Lemma \ref{lem4} that $h: V \rightarrow \sM_{g,n}^{ame}$ is regular in a neighbourhood of $p$. On the other hand, by Corollary \ref{cor6} we also have the natural regular map $\sM_{g,n}^{ame} \rightarrow \oM_{g,n}$, so it follows we have an isomorphism in a neighbourhood of $p$. Since $p$ was chosen arbitrarily we get $\sM_{g,n}^{ame} \cong \oM_{g,n}$.

\subsection{Extending Families of Stable Curves}

We can tweak Theorem \ref{thm6} to give necessary and sufficient criteria for being able to extend families of stable curves -- i.e. being able to extend a morphism $U \rightarrow \sM_{g,n}$ to a regular map $S \rightarrow \soM_{g,n}$. 

\begin{prop}\label{prop9} A family of smooth, pointed curves over $U \subset S$ extends to a family of stable, pointed curves over a neighbourhood of $p \in S \setminus U$ if and only if the local monodromy around $p$ is abelian and generated by Dehn twists (about multi-loops).
\end{prop}
\begin{proof}
Suppose the family extends over a neighbourhood of $p$ to a family of stable curves with the fibre over $p$ some pointed curve $C$. Denote by $V$ the local versal deformation space of $C$. The locus in $V$ corresponding to singular curves is a divisor $E$ with simple normal crossing at $[C] \in V$. Each component of $E$ corresponds to a node of $C$ so that the monodromy associated to it is a Dehn twist about the vanishing cycle corresponding to the node. This means that $\pi_1(V \setminus E) \cong \Z^k$ (where $k$ is the number of nodes in $C$) is generated by Dehn twists about disjoint loops. In particular, every element of $\pi_1(V \setminus E)$ is a Dehn twist about some multi-loop. Hence, the local monodromy around $p$, which is a subgroup of $\pi_1(V \setminus E)$, is abelian and generated by Dehn twists. 

Conversely, if the local monodromy around $p$ is abelian and generated by Dehn twists then by Theorem \ref{thm6} $U \rightarrow \sM_{g,n}$ extends to a regular map $S \rightarrow \oM_{g,n}$ in a neighbourhood of $p$. We want to show that we can lift $S \rightarrow \oM_{g,n}$ to $S \rightarrow \soM_{g,n}$. 

Denote by $C$ the curve corresponding to the image of $p$ in $\oM_{g,n}$ and let $W$ be the local versal deformation space of $C$. The neighbourhood of $[C] \in \soM_{g,n}$ is isomorphic to the quotient stack $[W/\Aut(C)]$ where $\Aut(C)$ is the automorphism group of $C$. After possibly restricting $S$ the problem comes down to lifting $S \rightarrow W/\Aut(C)$ to a morphism $S \rightarrow W$. 

If we denote by $G$ the fundamental group of the intersection of $\sM_{g,n}$ and a small neighbourhood of $[C] \in \soM_{g,n}$ then we have the short exact sequence
$$1 \rightarrow \Z^k \rightarrow G \xrightarrow{\phi} \Aut(C) \rightarrow 1$$
where $k$ is the number of nodes in $C$. In order to be able to lift $S \rightarrow W/\Aut(C)$ to $W$, the monodromy around $p$ (which is the image of $\pi_1(V \cap U) \rightarrow G \rightarrow \pi_1(\sM_{g,n})$ where $V$ is a small neighbourhood of $p$) must lie in the kernel of $\phi$. But the kernel of $\phi$ is precisely the subgroup of $G$ generated by Dehn twists around the $k$ vanishing cycles corresponding to the nodes of $C$. So if the monodromy is generated by Dehn twists it lies in $\mbox{ker}(\phi)$ and we are done.
\end{proof}

\begin{cor}\label{cor7} Let $S$ be a normal variety and $U \subset S$ a dense open subvariety. A morphism $U \rightarrow \sM_{g,n}$ extends in a neighbourhood of $p \in S \setminus U$ to a regular map $S \rightarrow \oM_{g,n}$ if and only if the local monodromy around $p$ is virtually abelian. 
\end{cor}
\begin{proof}
One direction is implied by Theorem \ref{thm6} so we just need to consider what happens if $U \rightarrow \sM_{g,n}$ extends around $p$. In this case there exists a finite cover $f: \tS \rightarrow S$ such that $S \rightarrow \oM_{g,n}$ lifts to a morphism $\tS \rightarrow \soM_{g,n}$. This follows from the same argument as in the proof of Theorem \ref{thm6} in section \ref{ss4}. Namely, $\soM_{g,n}$ has a finite flat cover $Z \rightarrow \soM_{g,n}$ by a projective scheme $Z$ (see \cite{kv}) and we take $\tS = S \times_{M_{g,n}} Z$. 

Then, by Proposition \ref{prop9}, the local monodromy around any point of $f^{-1}(p)$ is abelian so the local monodromy around $p$ must have been virtually abelian. 
\end{proof}

\begin{Remark}\label{rem:jo}
As a consequence to Corollary \ref{cor7} we get the following result first proved in \cite{jo}:

\begin{cor}\label{cor4} Let $D = S \setminus U$ be a normal crossing divisor at $p$. The morphism $U \rightarrow \sM_{g,n}$ extends in a Zariski neighbourhood of $p$ to a regular map $S \rightarrow \soM_{g,n}$ if and only if it extends over the generic points of $D$.
\end{cor}
\begin{proof}
If $D = S \setminus U$ is a normal crossing divisor at $p$ then in a small neighbourhood of $p$ the fundamental group of $U$ is free abelian. Thus $U \rightarrow \sM_{g,n}$ extends to $S \rightarrow \oM_{g,n}$. 

Now look at $\soM_{g,n} \rightarrow \oM_{g,n}$ and consider $f: S \times_{\oM_{g,n}} \soM_{g,n} \rightarrow S$. Locally in the \'etale topology this map looks like $f: [Y/G] \rightarrow Y/G = S$ for some scheme $Y$. If $Y$ is not normal then we replace it by its normalization $\tY$. The action of $G$ on $Y$ lifts to an action of $G$ on $\tY$ and we instead consider the map $f: [\tY/G] \rightarrow \tY/G$. Notice that $\tY/G \rightarrow Y/G = S$ is a birational, quasi-finite map and $S$ is smooth so that we get $\tY/G = S$. We would like to show that $f$ is an isomorphism.  

The locus where $f: [\tY/G] \rightarrow \tY/G = S$ is an isomorphism is where $g: \tY \rightarrow \tY/G = S$ is \'etale. Now $S$ is smooth, $\tY$ is normal and $g$ is \'etale over an open subset of $S$ whose complement has codimension $\ge 2$ (this is where we use that the morphism $U \rightarrow \sM_{g,n}$ extends over the generic points of $D$). By the purity of branch locus theorem (see for example p.461 of \cite{ak}) this implies $g$ is \'etale everywhere. Hence $f$ is an isomorphism. 

Thus we get a regular map $S = [\tY/G] \rightarrow \soM_{g,n}$ which lifts $S \rightarrow \oM_{g,n}$. 

\end{proof}

For much the same reason, we also get the extension theorem from \cite{mo} which generalized \cite{jo} to the case when the base is log regular (i.e. schemes with boundaries like those of toric varieties). In a slightly different direction, see also \cite{sa} for a discussion about extending to a log-smooth family instead of a stable family. 
\end{Remark}

\begin{Remark}
In \cite{st} Stix considers the analogous problem of extending $U \rightarrow \sM_g$ to a morphism $S \rightarrow \sM_g$ (no compactifying). Similarly, in \cite{bo} Boggi considers extending $U \rightarrow \sM_g$ to a morphism $S \rightarrow \tilde{\sM}_g$ where $\tilde{\sM}_g$ is the partial compactification consisting of stable curves of compact type. In both of these cases the condition for an extension to exist can be expressed in terms of monodromy. 
\end{Remark}

\section{The Moduli Space of Principally Polarized Abelian Varieties}\label{s4}

Denote by $\sA_g$ the moduli stack of $g$-dimensional principally polarized abelian varieties and by $A_g^{BB}$ its Baily-Borel (Satake) compactification. In this section we prove: 

\begin{thm}\label{thm5}
$A_g^{BB}$ is the maximal projective AME compactification of $\sA_g$. 
\end{thm}

\begin{Remark}
In the case of $\sA_2$ we can show that $A_2^{BB}$ is actually the maximal AME compactification (not just the maximal projective AME compactification). It is probably true more generally that $\sA_g^{ame} = A_g^{BB}$ but we cannot prove this at the moment. 
\end{Remark}

\subsection{Degeneration of Hodge Structures}

We follow the discussion from \cite{cat}. A (principally) polarized Hodge structure of weight one is a $\Z$ lattice $H_\Z$ of rank $2g$ for some $g \ge 0$ whose complexification $H = H_\Z \otimes \C$ is equipped with a decreasing Hodge filtration $H = F^0 \supset F^1 \supset F^2 = 0$ such that $H = F^1 \oplus \overline{F^1}$ ($F^1$ is often denoted $H^{1,0}$). The (principal) polarization is given by a (unimodular) non-degenerate, skew-symmetric bilinear form $\ip{\cdot}{\cdot}: H_\Z \times H_\Z \rightarrow \Z$ satisfying $\ip{F^1}{F^1}=0$ and $i\ip{v}{\ov} > 0$ for any non-zero $v \in F^1$. 

Fixing $H$ and $\ip{\cdot}{\cdot}$, the classifying space $D$ of polarized Hodge structures is the open subset of 
$$\check{D} = \{F^1 \in G(g,H): \ip{F^1}{F^1}=0\} \subset G(g,H)$$
subject to the condition $i\ip{v}{\ov}>0$ for any non-zero $v \in F^1$ (here $G(g,H)$ denotes the Grassmannian of $g$-planes in $H$). One can realize $D$ as the Siegel upper half space 
$$\H_g = \{Z_{g \times g}: Z=Z^t, \mbox{Im}(Z) > 0\}$$
as follows. Fix an integral symplectic basis $e_1, \dots, e_g, f_1, \dots, f_g$ of $H$ with respect to $\ip{\cdot}{\cdot}$. Because $i \ip{v}{\ov}>0$ for any $v \in F^1$ one can find a basis $\{v_1, \dots, v_g\}$ of $F^1$ such that $v_j = f_j + \sum_k Z_{kj} e_k$ where $\mbox{Im}(Z) > 0$. The condition $\ip{F^1}{F^1}=0$ implies $Z=Z^t$. $Z$ is called the normalized period matrix. 

Denote by $\Delta := \{z \in \C: |z| < 1\}$ the open unit disk and by $\Delta^* := \Delta \setminus \{0\}$ the punctured unit disk. Consider a variation of weight one polarized Hodge structures over $\Delta^*$. Taking a point to its normalized period matrix gives a period map $\phi: \Delta^* \rightarrow \Gamma \diagdown D$ where 
$$\Gamma \cong Aut(H_\Z, \ip{\cdot}{\cdot}) = \{\sigma: H_\Z \rightarrow H_\Z: \ip{v}{w} = \ip{\sigma v}{\sigma w} \text{ for any } v,w \in H_\Z \} \cong Sp_{2g}(\Z).$$
Denote by $\widetilde{\Delta^*}$ the universal cover of $\Delta^*$. If we identify $\widetilde{\Delta^*}$ with the upper half plane then the covering map $\widetilde{\Delta^*} \cong \H_1 \rightarrow \Delta^*$ is given by $z \mapsto \exp(2\pi iz) =: t \in \Delta^*$. Let $\widetilde{\phi}: \widetilde{\Delta^*} \rightarrow D$ be a lift of $\phi$ to the universal cover of $\Delta^*$. Denote by $\widetilde{Z}(z)$ and $Z(t)$ the normalized period matrices of $\widetilde{\phi}(z)$ and $\phi(t)$ (note that $\phi(t): \Delta^* \rightarrow \Gamma \diagdown D$ so one can think of $\phi(t): \Delta^* \rightarrow D$ and $Z(t)$ as multivalued functions). 

Denote by $T \in \Gamma$ the image of $1$ under the map $\phi_*: \Z \cong \pi_1(\Delta^*) \rightarrow \pi_1(\Gamma \diagdown D) \cong \Gamma$. By Landman's monodromy theorem \cite{lan}, $(T^k-I)^2=0$ for some $k \in \Z$. So after a finite base change we can assume $(T-I)^2=0$. Let $N = log(T) = T-I$. Since $N^2=0$ we have $\mbox{Im}(N) \subset \mbox{ker}(N)$. Then the monodromy weight filtration $W(N)$ induced by $N$ is 
$$0 \subset W_0 \subset W_1 \subset W_2 = H$$
where $W_0 = \mbox{Im}(N)$ and $W_1 = \mbox{ker}(N)$. An elementary computation shows $\ip{v}{Nw} = \ip{-T^{-1}Nv}{w}$ so, since $\ip{\cdot}{\cdot}$ is non-degenerate, $\ip{v}{Nw}=0$ for all $w \in H$ if and only if $Nv=0$. This means $W_0^\perp = W_1 \supset W_0$ and so $W_0$ is isotropic. 

Choose the symplectic basis above such that $W_0 = \span \{e_1, \dots, e_r\}$ where $r = \mbox{rk}(N)$. In this basis 
$N = \left( \begin{matrix} 0 & \eta \\ 0 & 0 \end{matrix} \right)$ and
$\eta = \left( \begin{matrix} \eta' & 0 \\ 0 & 0 \end{matrix} \right)$ 
where $\eta'$ is an $r \times r$ symmetric (since $T \in Sp_{2g}(\Z)$) matrix. Define $\widetilde{\psi}: \H_1 \rightarrow \check{D}$ by $\widetilde{\psi}(z) = \exp(-zN) \widetilde{\phi}(z)$. Since $\widetilde{\phi}(z+1)=T \widetilde{\phi}(z)$ we find $\widetilde{\psi}(z+1)=\widetilde{\psi}(z)$ so $\widetilde{\psi}$ descends to a map $\psi: \Delta^* \rightarrow \check{D}$. Denote by $\widetilde{W}(z)$ and $W(t)$ the normalized period matrices of $\widetilde{\psi}(z)$ and $\psi(t)$. Since $\exp(-zN) = \left( \begin{matrix} I & -z\eta \\ 0 & I \end{matrix} \right)$ we get $\widetilde{W}(z) = \widetilde{Z}(z) - z \eta$ or equivalently $Z(t) = W(t) + \frac{1}{2\pi i}(\log t) \eta$. 

The nilpotent orbit theorem \cite{sch} shows that $\psi: \Delta^* \rightarrow \check{D}$ has a removable singularity at zero which means $W(t)$ is holomorphic on $\Delta$. Write $W(t) = \left( \begin{matrix} W_{11}(t) & W_{12}(t) \\ W_{21}(t) & W_{22}(t) \end{matrix} \right)$ where $W_{22}(t)$ is a $(g-r) \times (g-r)$ symmetric matrix. 

\begin{lem} \label{lem8} Let $\Gr^W_1 = W_1/W_0$. The weight one (principally) polarized Hodge structure on $H$ induces a weight one (principally) polarized Hodge structure on $\Gr^W_1$ with normalized period matrix $W_{22}(t)$. 
\end{lem}
\begin{proof}
$W_1 = \span \{e_1, \dots, e_g, f_{r+1}, \dots, f_g\}$ while $W_0 = \span \{e_1, \dots, e_r\}$ which means 
$$\Gr^W_1 = \span \{[e_{r+1}], \dots, [e_g], [f_{r+1}], \dots, [f_g]\}.$$ 
Since $F^1 = \span \{v_j = f_j + \sum_k Z_{kj}(t) e_k: j = 1, \dots, g \}$ we find $F^1 \cap W_0 = \emptyset$ and $F^1 \cap W_1 = \span \{v_{r+1}, \dots, v_g \}$. Using the fact $\mbox{Im}(Z) > 0$ it is then elementary to see 
$$\Gr^W_1 = [F^1 \cap W_1] \oplus \overline{[F^1 \cap W_1]}$$
so that $F^1 \cap W_1$ defines a weight one Hodge filtration on $\Gr^W_1$. 

Since $W_0$ is isotropic $\ip{\cdot}{\cdot}$ descends to give a (principal) polarization on $W_0^\perp/W_0 = \Gr^W_1$. Fixing the symplectic basis $[e_{r+1}], \dots, [e_g], [f_{r+1}], \dots, [f_g]$ for $\Gr^W_1$, 
$$[F^1 \cap W_1] = \span \{ [v_j] = [f_j] + \sum_{k > r} Z_{kj}(t) [e_j]: j=r+1, \dots, g\}.$$
So, if we write $Z(t) = \left( \begin{matrix} Z_{11}(t) & Z_{12}(t) \\ Z_{21}(t) & Z_{22}(t) \end{matrix} \right)$ where $Z_{22}(t)$ is a $(g-r) \times (g-r)$ symmetric matrix, the normalized period matrix for $\Gr^W_1$ is $Z_{22}: \Delta^* \rightarrow D$. But $Z_{22}(t)$ extends to a holomorphic map $\Delta \rightarrow \check{D}$ because $Z_{22}(t) = W_{22}(t)$. 

A priori $Z_{22}(0) \in \overline{\H}_{g-r}$ where $\overline{\H}_g = \{Z_{g \times g}: Z=z^t, \mbox{Im}(Z) \ge 0\}$. However, for every non-zero real $(g-r)$ column matrix $X$ we have that $X^t Z_{22}(t) X$ lies in $\overline{\H}_1$ if $t \in \Delta$ and in $\H_1$ if $t \in \Delta^*$. So by the open mapping theorem $X^t Z_{22}(t) X \in \H_1$ for any $t \in \Delta$. Since this is true for any $X$ as above we get $Z_{22}(0) \in \H_{g-r}$. 
\end{proof}

\subsection{Proof of Theorem \ref{thm5}: $\sA_g^{amep} = A_g^{BB}$}

\subsubsection{Part I: Proof that $A_g^{BB}$ has AME property}

By Lemma \ref{lem7}, to show $(\sA_g,A_g^{BB})$ has the AME property it suffices to consider normal surfaces $U \subset S$. Suppose $S$ is a small neighbourhood of $p \in D = S \setminus U$ such that the local monodromy around $p$ with respect to a given morphism $f|_U: U \rightarrow \sA_g$ is virtually abelian. We need to show $f|_U$ extends to a regular map $f: S \rightarrow A_g^{BB}$. Lemma \ref{lem4} allows us to replace $U$ be any finite cover so that we can assume the local monodromy around $p$ is abelian. 

Since $S$ is normal it has isolated singularities so, after restricting $S$ to a smaller neighbourhood if necessary, one can assume $S$ is smooth with $D$ a normal crossing divisor except possibly at $p$. Let $\pi: S' \rightarrow S$ be a log-resolution of the pair $(S,D)$ -- meaning $\pi$ is a proper, birational morphism such that $S'$ is smooth and $D' = \pi^{-1}(D) \cup \{\mbox{exceptional divisors}\}$ is a normal crossing divisor. By blowing up further if necessary one can assume $S' \setminus D' = U' \xrightarrow{f|_U \circ \pi} \sA_g$ extends to a regular map $f': S' \rightarrow A_g^{BB}$. 

Let $C$ be a component in the exceptional divisor above $p$ and denote by $q_1, \dots, q_m$ the intersection points of $C$ with all the other divisors in $D'$. Denote by $U_C$ the intersection of $U'$ with a small open neighbourhood of $C$. Notice the image of the monodromy map $f'_*: \pi_1(U_C) \rightarrow \pi_1(\sA_g)$ is abelian. We will use this fact to show that $f'$ collapses $C$ to a point. This will mean that $f'$ collapses the whole exceptional divisor of $\pi$ to a point and so $f: S \rightarrow A_g^{BB}$ must be regular. 

The first cohomology group of the family of principally polarized abelian varieties corresponding to $U_C \rightarrow \sA_g$ gives a variation of principally polarized Hodge structures of weight one, which we denote $H$, on $U_C$. Let $Z$ be the normalized period matrix of $H$. 

Since $D'$ is normal crossing the kernel of the map $j: \pi_1(U_C) \rightarrow \pi_1(V_C)$ where $V_C = U_C \cup C \setminus \{q_i\}$ is isomorphic to $\Z$ and generated by a small loop around $C$. Denote by $T$ the image of this generator under the monodromy map $\pi_1(U_C) \rightarrow \pi_1(\sA_g)$ and let $N = T-I$. Notice that $T$ is only defined up to conjugacy since we can move the base point of $U_C$. However, because the image of the monodromy map is abelian we get a well defined action $N: H \rightarrow H$. Let $W^0 = \mbox{Im}(N)$ and $W^1 = \mbox{ker}(N)$. After taking an appropriate cover branched along $C$, Landman's monodromy theorem says $N^2=0$ so we have a filtration $0 \subset W^0 \subset W^1 \subset H$ and a quotient $\Gr = W^1/W^0$. 

Locally around any point $q \in C \setminus \cup_i q_i$ of $C$ choose coordinates $t_1,t_2$ such that $C$ is the locus $t_1=0$. By Lemma \ref{lem8}, the normalized period matrix of $H$ is $Z(t_1,t_2) = W(t_1,t_2) + \frac{1}{2\pi i} (\log t_1) \eta$ where $W(t_1,t_2) \in \H_g$ is holomorphic. Then  $H$ induces a weight one principally polarized Hodge structure on $\Gr = W^1/W^0$ with normalized period matrix $W_{22}(t_1,t_2) \in \H_{g-r}$ where $r = \mbox{rk}(N)$. 

Notice that $\Gr$ together with its Hodge structure is well defined over $V_C$ so we get a morphism $g: V_C \rightarrow \sA_{g-r}$. Since the monodromy of $H$ over $U_C$ is abelian the induced monodromy on $\Gr$ over $U_C$ is also abelian. Thus the image of $g_*: \pi_1(V_C) \rightarrow \pi_1(\sA_{g-r})$ is abelian. Restricting to $C \setminus \cup_i q_i \subset V_C$ we get a weight one principally polarized Hodge structure $\Gr$ whose associated period map $g_*: \pi_1(C \setminus \cup_i q_i) \rightarrow \pi_1(\sA_{g-r}) = Sp_{2(g-r)}(\Z)$ has abelian image.  But by Corollary 4.2.9 of \cite{del}, the connected component containing the identity of the Zariski closure of this image is semi-simple. Arguing as in the proof of \ref{prop8} shows that the monodromy on $C \setminus \cup_i q_i$ is finite. Hence, the morphism $g|_C: C \setminus \cup_i q_i \rightarrow \sA_{g-r}$ collapses $C$ to a point. But, by the construction of $A_g^{BB}$, the morphism $f': S' \rightarrow A_g^{BB}$ maps $C \setminus \cup_i q_i$ to the boundary component $A_{g-r} \subset A_g^{BB}$ via $g|_C$. So $f'$ collapses $C$ to a point. 

\subsubsection{Part II: Proof that $A_g^{BB}$ is maximal projective}

It remains to show that $A_g^{BB}$ is maximal among all projective AME compactifications of $\sA_g$. We make use of $A_g^{perf}$, the perfect (1st Voronoi) toroidal compactification of $A_g$ (see \cite{amrt}). By Lemma \ref{lem9} below we have a sequence of regular maps $A_g^{perf} \xrightarrow{\mu_1} A_g^{amep} \xrightarrow{\mu_2} A_g^{BB}$. The case $g=1$ is trivial since $\sA_1 \cong \sM_{1,1}$. We will first deal with the case $g=2$ and then use it to prove $A_g^{BB}$ is maximal for arbitrary $g$. 

Recall that $A_2^{BB} = A_2 \amalg A_1 \amalg A_0$ while $A_2^{perf}$ is isomorphic to $\oM_2$ via the Torelli map. The restriction of $\mu = \mu_2 \circ \mu_1$ over $A_1 \amalg A_0$ is the map $f: [\sM_{1,2}/(\Z/2\Z)] \cong \oM_{1,2} \rightarrow A_1^{BB}$ given by forgetting the second marked point. Let $E$ be any elliptic curve. We will now show that the map $\mu_1$ must collapse the fibre $f^{-1}([E])$ to a point. 

To do this we construct a rational map from a surface to $\sA_2$ with abelian local mondromy around a point $Q$ whose total transform inside $\oM_2$ contains the fibre $f^{-1}([E])$. The key is to find a morphism from a surface to $\oM_2$ whose image contains $f^{-1}([E])$ but where the local monodromy around the preimage of $f^{-1}([E])$ is abelian (then we blow down the preimage to get $Q$). We proceed as follows. 

\begin{figure}[ht]
\begin{center}
 \psfrag{pi}{$\pi$}
 \psfrag{E}{$E$}
 \psfrag{B}{$B$}
 \psfrag{S}{$S$}
 \psfrag{C}{$C$}
 \psfrag{P}{$P$}
 \psfrag{F}{$F$}
 \hspace{0 cm}\puteps[0.50]{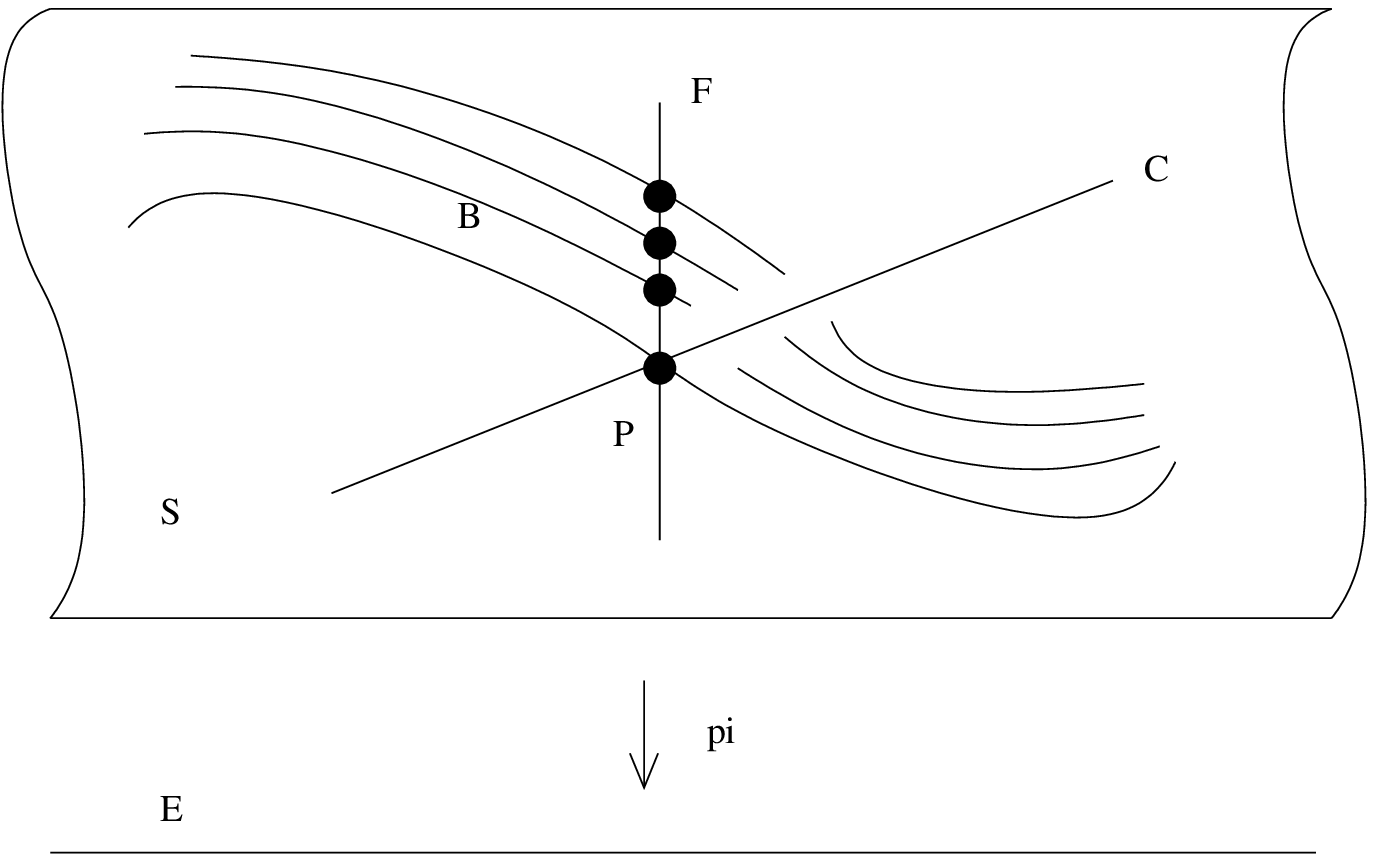}
\end{center}
\caption{}\label{f1F}
 \end{figure}

Fix a point $p \in E$ and denote by $\pi_1: E \times E \rightarrow E$ the projection onto the first factor. Let 
$V = \pi_{1*}(\O_{E \times E}(E \times p + \Delta))$ 
where $E \times p = \{(q,p): q \in E\} \subset E \times E$ and $\Delta$ is the diagonal. Alternatively, $V$ is a rank two vector bundle defined as the unique non-trivial extension $\O_E \rightarrow V \rightarrow  \O_E(p)$. Let $S$ be the surface $\P(V^\vee)$ and $\pi: S \rightarrow E$ the natural projection. The line bundle $\O_{E \times E}(E \times p + \Delta)$ defines a $2:1$ map $E \times E \rightarrow S$ mapping fibres to fibres. Denote by $B \subset S$ the branching locus of this map and by $C \subset S$ the image of $\Delta$ (or equivalently of $E \times p$). Notice that $\pi|_{B}: B \rightarrow E$ is an unramified, degree 4 map while $C$ is a section of $\pi$. Also, $B$ and $C$ meet only at (the missing) point $P \in S$ which is the image of $(p,p) \in E \times E$. Finally, denote by $F$ the fibre $\pi^{-1}\pi(P)$. A helpful picture to keep in mind is Figure \ref{f1F}. 

Let $U := S \setminus \{B \cup C \cup F\}$. To a point $q \in U$ we can associate the genus two curve obtained as the double cover of $\pi^{-1}\pi(q) \cong \P^1$ branched over $q$ and the five points $\pi^{-1}\pi(q) \cap (B \cup C)$ (by construction these six points are distinct). To do this globally over $U$ we consider the $\P^1$ bundle $\P(\pi^* V^\vee) \rightarrow U$ and take its double cover branched over the divisor corresponding to $\pi^{-1}\pi(q) \cap (B \cup C)$. 

To see why this is possible notice that we need to find a line bundle $L$ on $\P(\pi^* V^\vee)$ such that $L^{\otimes 2}$ is isomorphic to the branching divisor. Now 
$$\Pic(\P(\pi^* V^\vee)) \cong \Pic(U) \times \Z$$
and, under this isomorphism, the branching divisor has class $(\alpha, 6)$ for some $\alpha \in \Pic(U)$. Thus to find such a line bundle $L$ we just need to find a square root of $\alpha$. Now the pullback map $\Pic(E \setminus \{p\}) \rightarrow \Pic(U)$ is surjective since the fibres of $U \rightarrow E \setminus \{p\}$ are all isomorphic to $\P^1 \setminus \{ 5 \text{ points}\}$. Also, $\Pic(E \setminus \{p\})$ is a divisible group. Thus we can find a square root of $\alpha$. 

So now we have a family of hyperelliptic genus two curves $\Cu_U \rightarrow U$. The complement $D := S \setminus U$ is normal crossing except at $P$. Consider the log-resolution $S'$ and denote by $D'$ the preimage of $D$ which is now normal crossing (it's necessary to blow up at least twice since $B$ and $C$ are tangent at $P$ as one can check $B \cdot C = 2$ using the push-pull formula). If we denote by $U'$ the preimage of $U$ then the induced morphism $U' \rightarrow \sM_2$ extends to a regular map $S' \rightarrow \oM_2$ (this can be ensured either by Theorem \ref{thm6} or by blowing up further). 

\begin{figure}[ht]
\begin{center}
 \psfrag{p1}{\footnotesize{$p_1$}}
 \psfrag{p2}{\footnotesize{$p_2$}}
 \psfrag{p3}{\footnotesize{$p_3$}}
 \psfrag{p4}{\footnotesize{$p_4$}}
 \psfrag{p5}{\footnotesize{$p_5$}}
 \psfrag{E}{\footnotesize{$E$}}
 \psfrag{g1}{\footnotesize{$\gamma_1$}}
 \psfrag{g2}{\footnotesize{$\gamma_2$}}
 \psfrag{g=0}{\footnotesize{$g=0$}}
 \psfrag{g=1}{\footnotesize{$g=1$}}
 \psfrag{(a)}{\small{(a)}}
 \psfrag{(b)}{\small{(b)}}
 \psfrag{(c)}{\small{(c)}}
 \hspace{0 cm}\puteps[0.65]{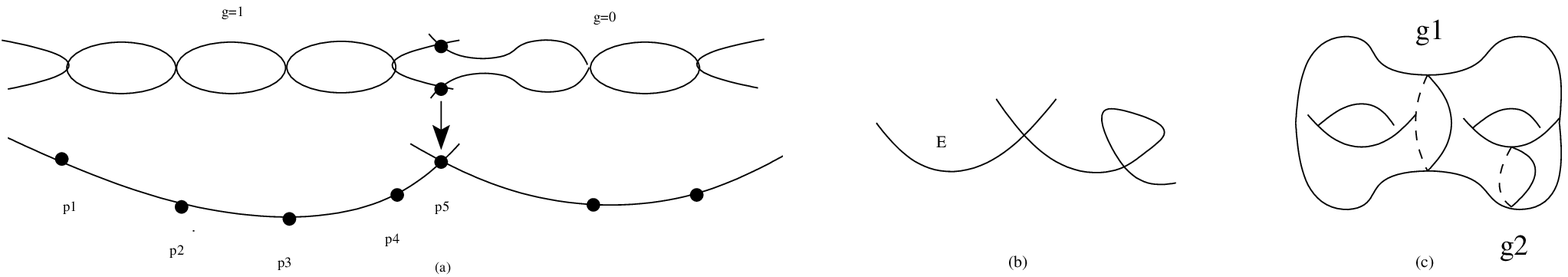}
\end{center}
\caption{}\label{f1E}
 \end{figure}

What is the image of $C' \subset S'$ (the proper transform of $C \subset S$) in $\oM_2$? If you approach a general point $c' \in C'$ along a general arc the associated genus two curves are double covers of $\P^1$ branched along six points, two of which are converging. The limiting admissible cover is illustrated in Figure \ref{f1E}(a). The points $p_1, \dots, p_4$ are independent of $c'$ (they are actually the branching points of $E \rightarrow \P^1$) and only $p_5$ varies as you move $c' \in C'$. This shows that $C'$ maps onto the fibre of $\oM_2 \rightarrow A_2^{BB}$ over $[E] \in A_1 \subset A_2^{BB}$. 

Now $C'$ intersects the exceptional divisor at one point which we call $P'$. The image of $P'$ in $\oM_2$ corresponds to the curve shown in Figure \ref{f1E}(b). The local monodromy around $P'$ induced by the morphism $U' \rightarrow \sM_2$ is generated by the two Dehn twists about loops $\gamma_1$ and $\gamma_2$ in Figure \ref{f1E}(c). However, as a map $U' \rightarrow \sA_2$ the Dehn twist about $\gamma_1$ induces the identity and so the local monodromy is generated only by $\gamma_2$. Consequently, if we denote by $W$ a small neighbourhood of $C' \subset S'$ then the map 
$$\pi_1(W \setminus C' \cup \{ \mbox{exceptional divisor} \}) \rightarrow \pi_1(\sA_2)$$
factors through $\pi_1(W \setminus C' \cup \{ \mbox{exceptional divisor} \}) \rightarrow \pi_1(W \setminus C')$. 

Since $W$ is smooth and $C' \cong E$ is a genus one curve, $W \setminus C'$ is homotopic to an $S^1$ bundle over $S^1 \times S^1$ so that $\pi_1(W \setminus C') \cong \Z^{\oplus 3}$. Hence the local monodromy around $C'$ induced by $U' \rightarrow \sA_2$ is abelian. This is the family we've been looking for. 

By blowing up sufficiently many points on $C'$ we can assume its proper transform $C'' \subset S''$ has an ample conormal bundle. Then we can blow down $C''$ to a point $Q \in \tS$. The local monodromy around $C''$ is a subgroup of the local monodromy around $C'$ and thus is abelian. Consequently, the local monodromy around $Q$ is abelian and we get a regular map $\tS \rightarrow A_2^{amep}$. We end up with a commutative diagram 
\begin{equation*}
\begin{CD}
C'' @>{\subset}>> S'' @>>> \oM_2 \\
@. @VVV @VVV \\
@. \tS @>>> A_2^{amep}
\end{CD}
\end{equation*}
which implies that $\mu_1: \oM_2 \rightarrow A_2^{amep}$ must collapse to a point the image of $C''$ inside $\oM_2$. This image is the same as that of $C'$ which is the fibre of $\oM_2 \rightarrow A_2^{BB}$ over $[E] \in A_1 \subset A_2^{BB}$. 

Finally, consider the restrictions of $\mu_1$ and $\mu_2$ to the preimage of $A_1^{BB} \subset A_2^{BB}$. We get $\oM_{1,2} \xrightarrow{\mu_1} \mu_2^{-1}(A_1^{BB}) \xrightarrow{\mu_2} A_1^{BB}$ where $\mu_1$ collapses to a point the fibre $\mu^{-1}([E])$. Since $\oM_{1,2}$ is irreducible and $A_1^{BB}$ is a normal curve the morphism $\mu: \oM_{1,2} \rightarrow A_1^{BB}$ is flat hence, by Lemma \ref{lem10}, $\mu_2: A_2^{amep} \rightarrow A_2^{BB}$ is an isomorphism over $A_1^{BB} \subset A_2^{BB}$ (we use that $A_2^{BB}$ is normal, projective). This shows $A_2^{amep} \cong A_2^{BB}$ (in fact it shows that $A_2^{ame} \cong A_2^{BB} \cong A_2^{amep}$). 

We will now use that $A_2^{amep} \cong A_2^{BB}$ together with the fact that $A_g^{amep}$ is projective to conclude that $A_g^{amep} \cong A_g^{BB}$ for $g > 2$. Recall that we have $A_g^{BB} = A_g \amalg A_{g-1} \amalg A_{g-2}^{BB}$ while the morphism $\mu: A_g^{perf} \rightarrow A_g^{BB}$ restricted to the preimage of $A_{g-1} \subset A_g^{BB}$ is isomorphic (after taking a $\Z/2\Z$ quotient) to the universal family $X_{g-1} \rightarrow A_{g-1}$. The preimage $\mu^{-1}(A_g \amalg A_{g-1})$ is a partial compactification of $A_g$ present in all toroidal compactifications and commonly denoted $A_g^{part}$. 

Fix a principally polarized abelian variety $A'$ of dimension $g-2$ and consider the morphism $\sA_2 \rightarrow \sA_g$ given by $[A] \mapsto [A' \times A]$. This map extends to a regular map $A_2^{part} \rightarrow A_g^{part}$ whose restriction to $X_1 \rightarrow A_{g-1}$ is given by $[E,p] \mapsto [A' \times E, (a,p)]$ where $a \in A'$ is the identity. Consider now the commutative diagram 
\begin{equation*}
\begin{CD}
A_2^{part} @>>> A_g^{perf} \\
@VVV @VVV \\
A_2^{BB} \cong A_2^{ame} @>>> A_g^{amep}.
\end{CD}
\end{equation*}
Since the morphism $A_2^{part} \rightarrow A_2^{BB}$ contracts the fibres of $X_1 \rightarrow A_1 \subset A_2^{BB}$ it follows that the image of such a fibre in $A_g^{perf}$ must be contracted to a point by the map $\mu_1: A_g^{perf} \rightarrow A_g^{amep}$. So in the fibre $A' \times E$ over $[A' \times E] \in A_{g-1}$ the locus contracted by $\mu_1$ includes $\{[A' \times E, (a,p)]: p \in E\}$. This is a fibre of the projection $A' \times E \rightarrow A'$ so that by Lemma \ref{lem10} the whole fibre $A' \times E$ is collapsed to $A'$. This means that the dimension of the fibre of $\mu_2: A_g^{amep} \rightarrow A_g^{BB}$ over $[A' \times E] \in A_{g-1}$ is at most $g-2$. Thus $X_{g-1} \subset A_g^{part}$, which has dimension $g(g+1)/2-1$, is collapsed by $\mu_2$ to a variety of dimension at most $\dim(A_{g-1})+g-2 = g(g+1)/2-2$. In particular, the divisor $D = \overline{X_{g-1}} \subset A_g^{perf}$ is collapsed by $\mu_2$ to a subvariety of lower dimension. 

On the other hand, if $g \ge 2$ then $Pic_{\Q}(A_g^{perf}) \cong \Q(L) \oplus \Q(D)$ where $L$ is the determinant of the Hodge bundle (\cite{egh} p.397). This means that $Pic_{\Q}(A_g^{amep}) \cong \Q(L)$. Since $L$ is pulled back from $A_g^{BB}$ by $\mu_2$ this implies $A_g^{amep} \cong A_g^{BB}$ because otherwise $A_g^{amep}$ would fail to be projective. 

\begin{lem}\label{lem9} Let $(\sX,\oX)$ be an AME pair and $\oX_{tor}$ a toroidal compactification of $X$. Then there exists a regular map $\oX_{tor} \rightarrow \oX$ extending the identity map on $X$. 
\end{lem}
\begin{proof}
Locally, $X \subset \oX_{tor}$ is isomorphic to $(\C^*)^n \subset X_\sigma$ where $X_\sigma$ is normal affine (and translations by $(\C^*)^n$ extend to $X_\sigma$). This means that the local fundamental group around any point $p \in \oX_{tor} \setminus X$ is abelian and hence its image in $\pi_1(X)$ is abelian. Since $\oX$ is an AME compactification the morphism $\oX_{tor} \rightarrow \oX$ must be regular. 
\end{proof}

\begin{lem}\label{lem10} Consider morphisms $Y \xrightarrow{f} X' \xrightarrow{g} X$ between complete, normal varieties such that the composition $g \circ f$ is flat with connected fibres and $f$ is surjective. If $f$ collapses a single fibre $(g \circ f)^{-1}(p)$ to a point then $g$ is an isomorphism. 
\end{lem}
\begin{proof}
Since $f$ collapses $(g \circ f)^{-1}(p)$ to a point the fibre $g^{-1}(p)$ is a point (here we use that $f$ is surjective).  Thus the relative dimension of $g$ is zero. 

This means $\dim(X) = \dim(X')$. If any fibre $g^{-1}(q)$ is positive dimensional then the dimension of $(g \circ f)^{-1}(q)$ is greater than $\dim(Y) - \dim(X') = \dim(Y) - \dim(X)$. This contradicts the fact that $g \circ f$ is flat. Since $g \circ f$ has connected fibres $g$ must also have connected fibres. Also, since $g \circ f$ is flat $g$ must be surjective. So $g$ is surjective with connected fibres of dimension $0$. This means $g$ is quasi-finite and bijective. Since $g$ is also proper it is in fact finite. Also $g$ is bijective so it is birational. Thus $g$ is a finite birational map and $X$ is normal so $g$ is an isomorphism. 
\end{proof}

\begin{prob} What is $A_g^{ame}$?
\end{prob}

It seems very likely that $A_g^{ame} \cong A_g^{BB}$. To show this one could try to extend the argument above by studying more carefully the morphism $\mu: A_g^{perf} \rightarrow A_g^{BB}$ deeper into the boundary of $A_g^{BB}$. 

\begin{Remark}
There is no reason to consider only principally polarized abelian varieties. More generally, the Baily-Borel compactification of an arithmetic quotient of any bounded symmetric domain should be its maximal (projective) AME compactification. The fact that the Baily-Borel compactification is AME follows in much the same way as in the proof of \ref{thm5} except that you need to deal with higher weight Hodge structures. That it is maximal follows from the fact that the Baily-Borel compactification is the log canonical model (a result we hope to explain in future work). 
\end{Remark}

\section{Some final remarks}

\subsection{Why abelian groups?}

Our definition of the AME property involves the class $\sC$ of abelian groups. Why not consider another class $\sD$ of finitely generated groups (say solvable groups) and look for varieties which have the $\sD$ monodromy extension property?

This is actually a reasonable thing to try to do. One could, for instance, try to find pairs $(X, \oX)$ having the SME (the solvable monodromy extension) property. But this is another story. The story in this paper is that there are natural moduli spaces which satisfy the AME property. In the case of $\sM_{g,n}$ we even get an ``if and only if'' statement (Corollary \ref{cor7}). 

If one tried to find (say) $\sM_{g,n}^{sme}$ there are two possibilities. Either it does not exist or there exists a proper surjective morphism $\alpha: \oM_{g,n} \rightarrow \sM_{g,n}^{sme}$. In the latter case if $\alpha$ is an isomorphism then we get nothing new so one might as well stay with the AME property since it is easier to check than the SME property. If $\alpha$ is not an isomorphism then it would indeed be interesting to identify explicitly $\sM_{g,n}^{sme}$. However, as many would agree, $\oM_{g,n}$ is likely the most natural compactification of $M_{g,n}$ so it would still be favourable to consider $\sM_{g,n}^{ame}$ instead of $\sM_{g,n}^{sme}$ since the former is probably more natural. 

\subsection{Why only varieties?}

In this paper we restrict our attention to testing for the AME property using varieties $U \subset S$. Why not also consider more general (formal) schemes?

\subsubsection{Main reason}

In order to show that $(X, \oX)$ has the AME property one needs to check that any morphism $U \rightarrow X$ extends in a neighbourhood of a point $p$ whenever the local mondromy around $p$ is virtually abelian. This is easier to do if you only need to check pairs of varieties $U \subset S$ rather than pairs of more arbitrary schemes. Furthermore, as this paper shows it suffices to consider only pairs of varieties $U \subset S$ in order to get natural maximal AME compactifications such as $\sM_{g,n}^{ame} \cong \oM_{g,n}$. 

\subsubsection{Secondary reason}

Even if we keep the current definition of an AME compactification it may still be interesting to ask what happens if $U \subset S$ is a pair of more general schemes. One is tempted to consider either $S$ an integral, separated scheme which is not necessarily of finite type or a formal scheme. The main examples of such schemes would be the localization or the completion of a variety along a subvariety. The primary example of the latter is the algebraic disk $\Spf(\C[t])$ which is the completion of $\Spec(\C[t])$ along the ideal $(t)$. 

The two main questions to ask are:
\begin{enumerate}
\item How do you define the local monodromy for these more general schemes?
\item If $(X,\oX)$ has the AME property does the extension property still hold for these schemes?
\end{enumerate}

Let us sketch a possible answer to these questions. In the finite type case, another way of describing the monodromy around $p$ is to look at the closure of the image of $f: U \rightarrow \oX$ and to take a small neighbourhood $V$ of the proper transform of $p$. Then the local monodromy around $p$ is the image of $\pi_1(V \cap f(U)) \rightarrow \pi_1(X)$. 

The map $f: U \rightarrow X$ extends if and only if the proper transform of $p$ is a point and the proofs in this paper show that if the proper transform is not a point then the monodromy around it is not virtually abelian. So, as before, the map extends in a neighbourhood of $p$ if the monodromy is virtually abelian. Now, if $S$ is an integral, separated scheme but not necessarily of finite type then a variant of this definition should still work assuming that $\oX$ is of finite type. 

For example, suppose $S$ is the localization of $T'$ along $T \subset T'$ where $T,T'$ are varieties. Then a morphism $U \rightarrow X$ corresponds to a morphism $f: U' \rightarrow X$ from some open subset $U' \subset T'$. So the local monodromy around some $p \in T$ is defined as the image of $\pi_1(f(U') \cap V) \rightarrow \pi_1(X)$ where $V$ is a small neighbourhood of the proper transform of $p$. Then the same argument as above shows that $f$ extends to a regular map in a neighbourhood of $p$ if the monodromy is virtually abelian. 

A little more subtle is the case when $S$ is a formal scheme which is the completion of $T'$ along $T \subset T'$. Given a map $f: U \rightarrow X$ one should be able to find a subvariety $Y \subset X$ (at least locally) such that the image of $U$ is the completion of $Y$ along $f(T)$. Then we would define the local monodromy around a point $p \in T$ as the image of $\pi_1(V \cap Y) \rightarrow \pi_1(X)$ where $V$ is a local neighbourhood of the proper transform of $p$.

\subsection{What about characteristic $p>0$}

There is a counter example to Theorem \ref{thm:b} in positive characteristic given in \cite{cf}, page 192. Since the tame algebraic fundamental group of the complement of a normal crossing divisor is abelian (\cite{gm}), this means that Theorem \ref{specialthm2} cannot hold in positive characteristic if we change to using the tame algebraic fundamental group and consider formal neighbourhoods instead of analytic ones.

\end{document}